\documentclass[12pt,a4paper]{article}
\usepackage{amsmath}
\usepackage{amssymb}
\usepackage{amsthm}
\usepackage{amsfonts}
\usepackage{latexsym}

\theoremstyle{plain}
\newtheorem{thm}{Theorem}
\newtheorem{cor}{Corollary}[section]
\newtheorem{lem}{Lemma}[section]
\newtheorem{prop}{Proposition}[section]
\theoremstyle{definition}
\newtheorem{defn}{Definition}[section]

\theoremstyle{remark}

\title{Szeg\"{o} kernels, Toeplitz operators,\\and equivariant fixed point formulae}
\author{Roberto Paoletti\footnote{\noindent{\bf Address:}
Dipartimento di Matematica e Applicazioni, Universit\`a degli Studi
di Milano Bicocca, Via R. Cozzi 53, 20125 Milano,
Italy; {\bf e-mail}: roberto.paoletti@unimib.it }}
\date{}

\begin{document}

\maketitle

\section{Introduction}

The aim of this paper is to apply algebro-geometric Szeg\"{o}
kernels to the asymptotic study of a class of trace formulae in equivariant geometric
quantization and algebraic geometry.

Let $(M,J)$ be a connected complex projective manifold, of complex dimension d, and let $A$ be an
ample line bundle on it. Let, in addition,
$G$ be a compact and connected g-dimensional Lie group acting holomorphically on $M$, in such a way that the action can
be linearized to $A$. For every $k\in \mathbb{N}$, the spaces of global holomorphic sections
$H^0\left (M,A^{\otimes k}\right)$ are linear representations of $G$, and therefore may be
equivariantly decomposed over its irreducible representations.

More precisely, let $\mathfrak{g}$ be the Lie algebra of $G$, and
let $\Lambda\subseteq
\mathfrak{g}^*$ be a set of weights parametrizing the family of all finite-dimensional irreducible representations of $G$. For
every $\varpi\in \Lambda$, denote by $V_\varpi$ the corresponding $G$-module. Given $\varpi\in \Lambda$
and a linear representation
of $G$ on a finite dimensional vector space $W$, we shall denote by $W_\varpi\subseteq W$ the $\varpi$-isotype
of $W$, that is, the maximal invariant subspace of $W$ equivariantly isomorphic to a direct sum of
copies of $V_\varpi$.
For every $k\in \mathbb{N}$,
we then have equivariant direct sum decompositions
\begin{equation}
\label{eqn:unitary-direct-sum}
H^0\left (M,A^{\otimes k}\right)=\bigoplus _{\varpi \in \Lambda}H^0\left (M,A^{\otimes k}\right)_\varpi.
\end{equation}

We can find a $G$-invariant Hermitian metric $h$ on $A$ such that the unique
compatible connection has curvature $\Theta=- 2i\omega$, where $\omega$ is a K\"{a}hler form on $M$.
The choice of $\omega$ determines a volume form on $M$, with respect to which
$\mathrm{vol}(M)=\frac{\pi^{\mathrm{d}}}{\mathrm{d}!}\,\int _M c_1(A)^{\mathrm{d}}$.

By way of motivation, suppose that a reductive connected algebraic group $\widetilde{G}$
acts on a complex projective manifold $M$; then for
any line bundle $A$ on $M$ there exists a positive integer $l$ such that the action linearizes to $A^{\otimes l}$
\cite{mumf}, \cite{dol}. Let $G\subseteq \widetilde{G}$
be a maximal compact subgroup. If $A$ is ample and $h'$ is an Hermitian metric on it whose compatible connection
has K\"{a}hler normalized curvature $\omega'$ , we may replace $h'$ and $\omega'$ by their $G$-averages
$h$ and $\omega$. The action of $G$ on $M$ is then holomorphic and Hamiltonian with respect to $\omega$.

These choices induce natural $G$-invariant Hermitian structures on every space of global sections of
$A^{\otimes k}$,
and (\ref{eqn:unitary-direct-sum}) is a \textit{unitary} equivariant isomorphism.
For $k\in \mathbb{N}$ and $\varpi\in \Lambda$, we shall denote by
$$P_{\varpi,k}:H^0\left (M,A^{\otimes k}\right)\rightarrow H^0\left (M,A^{\otimes k}\right)_\varpi$$
the orthogonal projector. With abuse of language, if $\mathcal{C}^\infty\left (M,A^{\otimes k}\right)$
is the space of all smooth
global sections of $A^{\otimes k}$, we shall also denote by $P_{\varpi,k}$ the orthogonal projector
$\mathcal{C}^\infty\left (M,A^{\otimes k}\right)\rightarrow H^0\left (M,A^{\otimes k}\right)_\varpi$.

Let, furthermore, $\gamma:M\rightarrow M$ be a biholomorphism, also admitting a linearization to a
\textit{unitary} automorphism $\widetilde{\gamma}$ of $(A,h)$; in particular, $\gamma$ is a symplectomorphism
of $(M,\omega)$. We shall also denote by $\widetilde{\gamma}$ the induced linearization on
$A^{\otimes k}$, for every $k\in \mathbb{Z}$.
For every $k\in \mathbb{N}$, let
$$
\widetilde{\gamma}_k:H^0\left (M,A^{\otimes k}\right)\rightarrow H^0\left (M,A^{\otimes k}\right),
\,\,\,\,\,\,s\mapsto \widetilde{\gamma}\circ s\circ \gamma ^{-1}
$$
be the unitary automorphism induced by $\widetilde{\gamma}$.

We shall make the additional assumption that $\widetilde{\gamma}$ commutes
with the action of $G$ on $A$; this is equivalent to the condition that $\Phi$ be $\gamma$-invariant, and therefore
$\gamma\left (\Phi^{-1}(0)\right)= \Phi^{-1}(0)$.
Under these circumstances, $\widetilde{\gamma}_k$ preserves the decomposition (\ref{eqn:unitary-direct-sum});
in other words, for every $\varpi\in \Lambda$ we have
$$
\widetilde{\gamma}_k\left(H^0\left (M,A^{\otimes k}\right)_\varpi\right)
= H^0\left (M,A^{\otimes k}\right)_\varpi.
$$

As a simple example, consider the unitary action of $S^1$ on $\mathbb{C}^{\mathrm{d}+1}$
given by $t\cdot (z_0,\ldots,z_\mathrm{d})=:\left (t\,z_0,t^{-1}\,z_1,\ldots,t^{-1}\,z_\mathrm{d}\right)$,
and let $\Gamma\in U(\mathrm{d}+1)$ be any unitary diagonal matrix.
This induces an holomorphic Hamiltonian circle action on $\mathbb{P}^\mathrm{d}$,
commuting with the holomorphic symplectomorphism
$\gamma$ of $\mathbb{P}^\mathrm{d}$ induced by $\Gamma$; both have tautological
linearizations to the hyperplane line bundle.

More generally, assume given a holomorphic Hamiltonian action of a Lie group $H$ on $(M,J,\omega)$, linearizing to $A$;
for $h\in H$, let $\psi_h:M\rightarrow M$ and $\widetilde{\psi}_h:A\rightarrow A$ be the associated maps.
If $G\subseteq H$ is a compact and connected subgroup, and $h$ centralizes $G$, then the
restriction of the action to $G$ and $\gamma=:\psi_h$, $\widetilde{\gamma}=:\widetilde{\psi}_h$ satisfy the previous
hypothesis.

Let us now describe the main object of study of this paper:

\begin{defn}\label{defn:psi-k}
For $f\in \mathcal{C}^\infty(M)$, let $M_f:\mathcal{C}^\infty\left(M,A^{\otimes k}\right)\rightarrow
\mathcal{C}^\infty\left(M,A^{\otimes k}\right)$ be the multiplication operator $s\mapsto f\,s$.
Let $\gamma:M\rightarrow M$, $\widetilde{\gamma}:A\rightarrow A$ be as above.

\begin{enumerate}
  \item For every
$\varpi\in \Lambda$ and $k\in \mathbb{N}$, we introduce the equivariant Toeplitz operator
$$T_{f}^{(\varpi,k)}=:P_{\varpi,k}\circ M_f\circ P_{\varpi,k}:H^0\left (M,A^{\otimes k}\right)\rightarrow
H^0\left (M,A^{\otimes k}\right),$$
which we shall view as an endomorphism of $H^0\left (M,A^{\otimes k}\right)_\varpi$.
  \item More generally, we may consider for every
$\varpi\in \Lambda$ and $k\in \mathbb{N}$ the compositions
\begin{equation*}%\label{eqn:psi-k}
\Psi_{\varpi,k}=\Psi_{\varpi,k}(\gamma,f)=:\widetilde{\gamma}_k\circ T_{f}^{(\varpi,k)}:
H^0\left (M,A^{\otimes k}\right)_\varpi\rightarrow H^0\left (M,A^{\otimes k}\right)_\varpi.
\end{equation*}
\end{enumerate}

\end{defn}

We shall show that, under familiar assumptions in the theory of symplectic reductions,
the trace of $\Psi_{\varpi,k}$ admits an asymptotic expansion as $k\rightarrow +\infty$,
and explicitly describe its leading term.

Of course, in the action free case and with $f=1$ the Lefschetz fixed point formula of \cite{as}
gives an exact expression for $\mathrm{trace}(\widetilde{\gamma}_k)$, but even in this case it may be of some interest that
Szeg\"{o} kernels, and more precisely their scaling limits (\cite{bs}, \cite{bsz}, \cite{sz}),
provide a relatively elementary approach to the leading asymptotics.
%Also, in the action free case and with $\widetilde{\gamma}=\mathrm{id}_\mathbb{A}$, the limit of
%$k^{-\mathrm{d}}\,\mathrm{trace}\left(T_{f}^{(k)}\right)$ as $k\rightarrow +\infty$ has already been studied in \cite{ber},
%where actually $\omega$ is assumed to be non-degenerate, but perhaps not strictly positive.
On the other hand, the action free case with $\widetilde{\gamma}=\mathrm{id}$ has been studied in
\cite{bdmg}, and more recently, under wider hypothesis on the symplectic structure, in \cite{ber}.
If $f=1$ and $\widetilde{\gamma}$ is the identity, then the trace of $\Psi_{\varpi,k}$ computes
$\dim H^0\left(M,A^{\otimes k}\right)_\varpi$, a natural object of study in the
setting of symplectic reduction and geometric quantization since the landmarks \cite{gs-gq}, \cite{gs-hq};
in fact,
exact formulae for these dimensions are provided, for each given $k$, by the principle
$[Q,R]=0$ \cite{mein}. A rather elementary approach to the asymptotics for $k\rightarrow +\infty$
has
been given in \cite{p-mm}, based on microlocal techniques.

In the general case, where the linearization and the equivariant Toeplitz operator are considered on
the same footing, the leading coefficient in the asymptotic expansion for $\mathrm{trace}(\Psi_{\varpi,k})$ is the
product of the leading coefficient of the Lefschetz  fixed point formula on the symplectic reduction of $(M,\omega)$
and a certain $G$-average of $f$, with a weighting that depends on $\varpi$ and $\gamma$.

In order to state the result more precisely, we need to describe some invariants associated to $\gamma$, $\widetilde{\gamma}$
and the linearization of the $G$-action.

Given the linearization, the action of $G$ on $(M,2\omega)$ is Hamiltonian; let
$\Phi:M\rightarrow \mathfrak{g}^*$ be the corresponding moment map. We shall assume that $0\in \mathfrak{g}^*$
is a regular value of $\Phi$, and that $G$ acts freely on $\Phi^{-1}(0)$; with minor complications,
the arguments below apply however to the case where the (necessarily finite) stabilisers of the points in
$\Phi^{-1}(0)$ all have the same cardinality.

Let
\begin{equation}
\label{eqn:defn-di-M-0}
p:\Phi^{-1}(0)\longrightarrow M_0=:\Phi^{-1}(0)/G
\end{equation}
be the projection onto the symplectic reduction $M_0$ of $M$.
Thus $p$ is a principal $G$-bundle, and
the K\"{a}hler structure $(\omega,J)$
of $M$ descends in a natural manner
to the quotient K\"{a}hler structure
$(\omega_0,J_0)$ of $M_0$.

Let us list the definitions that will build up the statement of Theorem \ref{thm:main}.

\begin{defn}
Associated to $\gamma$ and the $G$-action we have the following objects.
\label{defn:componenti-luogo-fisso}
\begin{enumerate}
              \item Since it commutes with the $G$-action,
$\gamma$ descends to a holomorphic symplectomorphism of $(M_0,\omega_0,J_0)$,
that we shall denote by $\gamma _0:M_0\rightarrow M_0$.
              \item Let $F_1,\ldots,F_\ell\subseteq M_0$ be the connected components of its fixed locus,
$\mathrm{Fix}(\gamma_0)$: for every
$l=1,\ldots \ell$, $F_l$ is a complex submanifold of $M_0$ of, say, complex dimension
 $\mathrm{d}_l$. Let $\mathrm{c}_l=:(\mathrm{d}-\mathrm{g})-\mathrm{d}_l$
denote its complex codimension.
              \item If $l=1,\ldots,\ell$ and $r\in F_l$, let $N_{l,r}$ be the normal space to
$F_l\subseteq M_0$ at $r$, and let $\gamma_{r}:N_{l,r}\rightarrow N_{l,r}$
be the unitary map induced by the holomorphic differential of $\gamma_0$ at $r$.
Then $\mathrm{id}_{N_{l,r}}-\gamma_{r}^{-1}$ is non-singular, and
\begin{equation}
\label{eqn:determinant-factor}
c_l(\gamma)=:\det \!{}_{\mathbb{C}}\left (\mathrm{id}_{N_{l,r}}-\gamma_{r}^{-1}\right)
%=\overline{\det \!{}_{\mathbb{C}}\left(\mathrm{id}_{N_{l,r}}-\gamma_{r}\right)}
\end{equation}
is constant on $F_l$.
              \item If $m\in p^{-1}(F_l)$
for some $l$, since $G$ acts freely on $\Phi^{-1}(0)$ there exists a unique $g_m\in G$ such that
$\gamma(m)=\mu_{g_m}(m)$, where $\mu:G\times M\rightarrow M$ is the given action on $M$;
the conjugacy class of $g_m$ only depends on $l$.
\item If $l=1,\ldots,\ell$, $\chi_\varpi(g_m)$ does not depend on the choice of
$m\in p^{-1}(F_l)$; we shall set $\chi_\varpi(F_l)=:\chi_\varpi(g_m)$,
($m\in p^{-1}(F_l)$).

            \end{enumerate}

 \end{defn}

To see that the conjugacy class of $g_m$ only depends on $l$, set $\widetilde{F_l}=:p^{-1}(F_l)$
and suppose $m,n\in \widetilde{F_l}$. If $\pi(m)=
\pi(n)$, there exists $h\in G$ such that $n=\mu_h(m)$. Therefore,
$\mu_{h\,g_m\,h^{-1}}(n)=\mu_h\big(\gamma(m)\big)=\gamma\big(\mu_h(m)\big)=\gamma(n),$ whence
$g_{\mu_h(m)}=h\,g_m\,h^{-1}$, $\forall \,h\in G$, $m\in p^{-1}(F_l)$.
If $\pi(m)\neq \pi(n)$, let $\eta:[0,1]\rightarrow F_l$ be a smooth path such that
$\eta(0)=\pi(m)$, $\eta(1)=\pi(n)$.
The principal $G$-bundle $p$ has a natural connection; let $\eta^\sharp:
[0,1]\rightarrow \widetilde{F_l}$ be the unique horizontal lift of $\eta$ such that $\eta^\sharp(0)=m$.
Then $\mu_{g_m}\circ \eta^\sharp,\,\gamma\circ \eta^\sharp:[0,1]\rightarrow
\widetilde{F_l}$ are both horizontal lifts of $\eta$, and satisfy
$\mu_{g_m}\circ \eta^\sharp(0)=\gamma(m)=\gamma\circ \eta^\sharp(0)$; hence they are equal.
Therefore, $\mu_{g_m}\left(\eta^\sharp(1)\right)=\gamma\left( \eta^\sharp(1)\right)$, and so
$g_{\eta^\sharp(1)}=g_m$. Since $\eta^\sharp(1)$ is in the same obit as $n$,
we conclude by the previous considerations that  $g_m$ is conjugate to $g_n$.

\begin{defn}
\label{defn:f-mediata-su-G}
Given $f\in \mathcal{C}^\infty(M)$, we define
$\overline{f}\in \mathcal{C}^\infty\left (M_0\right)$ as the $G$-average
of $f$,
viewed as a smooth function on $M_0$. In other words,
$$
\overline{f}(m_0)=:\int _Gf\Big(\mu_g(m)\Big)\,d\nu(g)\,\,\,\,\,\,\,\,(m_0\in M_0),$$
where $m\in p^{-1}(m_0)\subseteq \Phi^{-1}(0)$.
\end{defn}

\begin{defn}
\label{defn:quotient-linearization}
The ample line bundle $A$ descends
to an ample line bundle $A_0$ on $M_0$, and the linearization $\widetilde{\gamma}$
descends to a linearization $\widetilde{\gamma}_0$
on $A_0$. If $l=1,\ldots,\ell$ there exists a unique $h_l\in S^1$ such that
$\widetilde{\gamma}_0(r):A_{0}(r)\rightarrow A_{0}(r)$
is multiplication by $h_l$ for every  $r\in F_l$,
where $A_{0}(r)$ is the fiber of $A_0$ at $r$.
\end{defn}

With the above notation, we then have:

\begin{thm}
\label{thm:main}
Suppose $\varpi\in \Lambda$,
$f\in \mathcal{C}^\infty(M)$ and $\gamma:M\rightarrow M$, with unitary linearization
$\widetilde{\gamma}:A\rightarrow A$, are given as above, so that $\widetilde{\gamma}$
commutes with the action of $G$ on $A$. Let $\Psi_{\varpi,k}$ be as in Definition \ref{defn:psi-k}. Then:
\begin{description}
  \item[i):] If $\Phi^{-1}(0)=\emptyset$, then
$
\Psi_{\varpi,k}=0$ for $k\gg 0$.
  \item[ii):] If
$\Phi^{-1}(0)\neq \emptyset$,
assume that $0\in \mathfrak{g}^*$ is a regular value of $\Phi$ and that $G$ acts freely on
$\Phi^{-1}(0)$. Then  as $k\rightarrow +\infty$ there is an asymptotic expansion
\begin{eqnarray*}
\lefteqn{\mathrm{trace}(\Psi_{\varpi,k})}\\
&\sim& \dim(V_\varpi)\,\sum_{l=1}^\ell \left (\frac k\pi\right)^{\mathrm{d}_l}\,\frac{h_l^k}{c_l(\gamma)}
\,\chi_\varpi(F_l)\,\int_{F_l}\overline{f}\,\mathrm{vol}_{F_l}\cdot\left (1+\sum _{a\ge 1}k^{-a/2}\,c_{\varpi la}\right).
\end{eqnarray*}
\end{description}

\end{thm}

As a test case, suppose that $f=1$ and $\widetilde{\gamma}$ is the identity map.
Then the asymptotic expansion of the Theorem reduces to:
$$
\dim H^0\left(M,A^{\otimes k}\right)_\varpi
\sim \dim(V_\varpi)^2\,\left(\frac k\pi\right)^{\mathrm{d}-\mathrm{g}}\,
\mathrm{vol}(M_0)\,\cdot\left (1+\sum _{a\ge 1}k^{-a/2}\,c_{\varpi la}\right).
$$
Up to a different normalization convention for the volume form, this agrees with
Theorem 2 of \cite{p-mm} (where only powers of $k^{-1}$ appear).

The proof is largely based on the microlocal theory of the Szeg\"{o} kernel in \cite{bs}, and
on its developments in \cite{z}, \cite{bsz},
\cite{sz}.

To motivate the role of asymptotic expansions for Szeg\"{o} kernels, let us
dwell again on the action free case, first with $\widetilde{\gamma}=\mathrm{id}$.
Thus, given $f\in \mathcal{C}^\infty(M)$, we are considering the asymptotics of the trace of the level-k Toeplitz operator $T_f^{(k)}=:P_k\circ M_f\circ P_k:H^0\left (M,A^{\otimes k}\right)
\rightarrow H^0\left (M,A^{\otimes k}\right)$; here $P_k$ is the level-k Szeg\"{o} projector, that is,
the full orthogonal projector onto $H^0\left (M,A^{\otimes k}\right)$.
An asymptotic expansion
in this case has been proved in \cite{bdmg}, \S 13. On the other hand, inserting the
diagonal asymptotic expansion for the level-k Szeg\"{o} kernel of \cite{z} in the Schwartz kernel of
$T_f^{(k)}$ leads to
\begin{equation}
\label{eqn:toeplitz-expansion}
\mathrm{trace}\left (T_f^{(k)}\right)\sim \left (\frac k\pi\right)^{\mathrm{d}}\,\int _Mf\cdot \mathrm{vol}_M+
\mathrm{L.O.T.}.
\end{equation}
Trying to adopt the same approach to the case of the Lefschetz fixed point formula of \cite{as}, thus now with $f=1$,
one is led however to consider the asymptotics of the Szeg\"{o} kernel over off-diagonal points in $M\times M$ of the form $\big(m,\gamma(m)\big)$,
and this motivates the appearance of scaling limits around the fixed locus of $\gamma$ into the picture.
In the general equivariant case, the asymptotic concentration of the equivariant Szeg\"{o} kernels determines a further localization around the zero locus of the moment map.

It is in order to conclude this introduction by emphasizing that there is a broader
scope for the methods and techniques appearing in this paper.
Firstly, although our present focus is on the simpler
holomorphic context, the following analysis could be generalized
to the symplectic almost complex category, in view of
the microlocal description in \cite{sz} of the almost complex
analogues of Szeg\"{o} kernels. On the other hand, the study of compositions akin to those in
Definition \ref{defn:psi-k} is particularly relevant to the theory of Toeplitz
quantization  \cite{z-id}, where one studies quantum maps associated to contact transformations,
whose underlying symplectic maps are generally not holomorphic.
In fact, all the ingredients of this article are already in place in the unitarization
process of \cite{z-id}, where the non-unitarity of quantum maps induced
by a contactomorphism is corrected
by composing with appropriate
Toeplitz operators, and in some cases certain trace formulae of this type are given.

\textbf{Acknowledgments.} I am very endebted to the referee for several stimulating comments on possible developments
and for suggesting various improvements to the exposition.

\section{Proof of Theorem \ref{thm:main}}

The first statement of the Theorem is an immediate consequence of the theory of \cite{gs-gq}. More precisely, since
$\Phi(M)\subseteq \mathfrak{g}^*$ is a compact subset, if $\Phi^{-1}(0)=\emptyset$ then
$\varpi\not\in k\,\Phi(M)$, hence $H^0\left(M,A^{\otimes k}\right)_{\varpi}=0$, for $k\gg 0$.

In order to prove the second statement, we shall lift the problem to the CR structure of the associated circle bundle.
Let $X\subseteq A^*$ be the unit circle bundle, with $S^1$-action  $r:S^1\times X\rightarrow X$. The
connection 1-form $\alpha\in \Omega^1(X)$ defines a contact structure and a volume form
on $X$. In terms of the latter, we shall tacitly identify (generalized)
densities, half-densities and functions
on $X$ and $X\times X$.

With this in mind, there is a natural Hermitian structure on $L^2(X)$, and for every $k\in \mathbb{N}$
there are standard unitary isomorphisms
$\mathcal{C}^\infty\left (M,A^{\otimes k}\right)\cong \mathcal{C}^\infty(X)_k$; the latter is the
space of all smooth functions on $X$ such that $f\big (r_t(x)\big)=t^k\,f(x)$, $\forall\,t\in S^1,\,x\in X$.
By restriction, we
obtain the unitary isomorphisms
$H^0\left (M,A^{\otimes k}\right)\cong H(X)_k$, where
$H(X)_k\subseteq \mathcal{C}^\infty(X)_k$ is the $k$-th isotypical component of the Hardy space of $X$.
Given $s\in \mathcal{C}^\infty\left (M,A^{\otimes k}\right)$, we shall denote its image by $\widehat{s}\in
\mathcal{C}^\infty(X)_k$.

Furthermore, given our assumptions, the action $\mu:G\times M\rightarrow M$ naturally lifts
to an action of $G$ on $X$,
$\mu_X:G\times X\rightarrow X$, as a group of contactomorphisms; to lighten notation, we shall often
write $\mu$ for $\mu_X$ where no misunderstanding seems likely.

Similarly, $\widetilde{\gamma}:A\rightarrow A$ induces a
contactomorphism $\gamma_X:X\rightarrow X$, and $\gamma_X\circ \mu_g=\mu_g\circ \gamma_X$, $\forall\,g\in G$.

Thus $G$ acts on $H(X)_k$ by pull-back, $g:f\mapsto f\circ \mu_{g^{-1}}$, and the isomorphisms
$H^0\left (M,A^{\otimes k}\right)\cong H(X)_k$ are equivariant for this action.
In terms of the equivariant unitary isomorphism $s\mapsto \widehat{s}$, we may rewrite (\ref{eqn:unitary-direct-sum})
as
\begin{equation}
\label{eqn:unitary-direct-sum-hardy}
H(X)_k=\bigoplus _{\varpi\in \Lambda}H(X)_{\varpi,k}.
\end{equation}

\noindent
Similarly, for all
$s\in \mathcal{C}^\infty\left (M,A^{\otimes k}\right)$
we have $\widehat{\widetilde{\gamma}_k(s)}=\widehat{s}\circ \gamma _X^{-1}\in \mathcal{C}^\infty(X)_k$.

If $\left\{s_j^{(\varpi,k)}\right\}_j$ is any orthonormal basis of $H_{\varpi,k}(X)$, then
$$\Pi_{\varpi,k}(x,y)=:\sum _j s_j^{(\varpi,k)}(x)\,\overline{s_j^{(\varpi,k)}(y)}\,\,\,\,\,
(x,y\in X)
$$
is the $(\varpi,k)$-equivariant Szeg\"{o} kernel, that is, the distributional kernel of the
orthogonal projector onto the subspace $H(X)_{\varpi,k}\subseteq H(X)$; in particular, it does not depend on
the choice of $\left\{s_j^{(\varpi,k)}\right\}_j$. Thus, since $\left\{s_j^{(\varpi,k)}\circ \gamma_X\right\}_j$
is also an orthonormal basis of $H_{\varpi,k}(X)$, we have
 $$\Pi_{\varpi,k}\big(\gamma _X(x),\gamma _X(y)\big)=\Pi_{\varpi,k}(x,y),\,\,\,\,\forall\,x,y\in X.$$

Let $\mathrm{dens}_{X\times X}$ and $\mathrm{dens}_{X}$  denote the volume densities of $X\times X$
and $X$, respectively. Then
\begin{eqnarray}\label{eqn:eqz-per-traccia}
%\lefteqn{}
\mathrm{trace}\left(\Psi_{\varpi,k}\right)&=&\int _{X\times X}\Pi_{\varpi,k}\left (\gamma_X^{-1}(x),y\right)\,f(y)\,\Pi_{\varpi,k}\left (y,x\right)
\,\mathrm{dens}_{X\times Y}(x,y),\nonumber\\
&=&\int _{X\times X}\Pi_{\varpi,k}\big (x,\gamma_X(y)\big)\,f(y)\,\Pi_{\varpi,k}\left (y,x\right)
\,\mathrm{dens}_{X\times X}(x,y)\nonumber \\
&=&\int _X \Pi_{\varpi,k}\left (\gamma_X^{-1}(y),y\right)\,f(y)\,\mathrm{dens}_{X}(y).
\end{eqnarray}

The strategy to determine the asymptotics of $\mathrm{trace}\left(\Psi_{\varpi,k}\right)$
is then to insert in (\ref{eqn:eqz-per-traccia}) the asymptotic expansion for the scaling limits of
$\Pi_{\varpi,k}$ determined in \cite{p}. To this end, we shall apply a number of reductions, at each
step disregarding a contribution to the integral in (\ref{eqn:eqz-per-traccia}) which is
$O\left (k^{-\infty}\right)$.

Let us define
$
R(\Phi)=:\big \{(m,n)\in M\times M:\,\Phi(m)=0,\,n\in G\cdot m\big\}
$, $I(\Phi)=:(\pi\times \pi)^{-1}\big(R(\Phi)\big)$.
In other words,
$$
I(\Phi)=\Big\{(x,y)\in X\times X:\Phi\circ \pi(x)=0,\,y\in (G\times S^1)\cdot x\Big\}.
$$

\begin{lem}\label{lem:inotorno-phi-0}
Uniformly on compact subsets of
$X\times X\setminus I(\Phi)$, as $k\rightarrow +\infty$ we have
$\Pi_{\varpi,k}(x,y)=O\left (k^{-\infty}\right)$.
\end{lem}

Lemma \ref{lem:inotorno-phi-0}, whose proof will be postponed, implies the following:
if we fix an arbitrarily small $G$-invariant tubular neighborhood
$V\subseteq M$ of
$\Phi^{-1}(0)$, perhaps after disregarding a rapidly decaying contribution we may replace the
integration over $X$ in (\ref{eqn:eqz-per-traccia}) by an integration over $\pi^{-1}(V)$.
We shall express this by writing
\begin{eqnarray}\label{eqn:eqz-per-traccia-v-times-v}
%\lefteqn{\mathrm{trace}\left(\Psi_{\varpi,k}\right)}
\mathrm{trace}\left(\Psi_{\varpi,k}\right)\sim\int _{\pi^{-1}(V)}\Pi_{\varpi,k}\left (\gamma_X^{-1}(y),y\right)\,f(y)\,\mathrm{dens}_{X}(y).\nonumber
\end{eqnarray}
In particular, we may assume without loss that $G$ acts freely on $\overline{V}$.

As a further reduction, let us define
\begin{eqnarray}
\label{eqn:defn-of-sk}
S_k&=:&\left\{m\in V:\mathrm{dist}_M\big(G\cdot m,G\cdot \gamma(m)\big)<
2\,k^{-2/5}\right\},\\
S_k'&=:&\left\{m\in V:\mathrm{dist}_M\big(G\cdot m,G\cdot \gamma(m)\big)>k^{-2/5}\right\}.
\nonumber
\end{eqnarray}
Let $\{\sigma_k,\sigma_k'\}$ be a partition of unity on $V$
subordinate to the open cover $\{S_k, S_k'\}$. We may assume
$\sigma_k(m)=\sigma\left (k^{2/5}\,\mathrm{dist}_M\big(G\cdot m,G\cdot \gamma(m)\big)\right)$, for a fixed
smooth function $\sigma:\mathbb{R}\rightarrow \mathbb{R}$. We shall write $\sigma_k$ for $\sigma_k\circ \pi$,
$\sigma_k'$ for $\sigma_k'\circ \pi$.
Inserting the equality $\sigma_k+\sigma_k'=1$ in
(\ref{eqn:eqz-per-traccia-v-times-v}), the integral splits in two summands.
One of these gives a negligible contribution to the asymptotics as $k\rightarrow +\infty$:

\begin{lem}
\label{lem:annullamento-primo-termine-s}
As $k\rightarrow +\infty$, we have
$$\int _{\pi^{-1}\left (S_k'\right)}\sigma_k'(y)\,\Pi_{\varpi,k}\left (\gamma_X^{-1}(y),y\right)\,f(y)
 \,\mathrm{dens}_{X}(y)=O\left (k^{-\infty}\right).$$
\end{lem}

\textit{Proof.} It suffices to show that $\left |\Pi_{\varpi,k}\left (\gamma_X^{-1}(y),y\right)\right|=O\left (k^{-\infty}\right)$
uniformly for $y\in \pi^{-1}\left (S_k'\right)$.
To this end, recall that
\begin{equation}
\Pi_{\varpi,k}\left (\gamma_X^{-1}(y),y\right)=
\dim(V_\varpi)\,\int _G\chi _\varpi\left (g\right)\,\Pi_k\left (\mu_{g}\circ
\gamma_X^{-1}(y),y\right)\,
d\nu(g).
\label{eqn:equivariant-projection}
\end{equation}
To simplify notation, let us write $\mathrm{dist}_M$ for the composition
$\mathrm{dist}_M\circ (\pi\times \pi):X\times X\rightarrow \mathbb{R}$.
If $y\in \pi^{-1}\left (S_k'\right)$,
then
$\mathrm{dist}_M\left (\mu_{g}\circ \gamma_X^{-1}(y),y\right)> k^{-2/5}$ for all $g\in G$.
By the off-diagonal estimates on the Szeg\"{o} kernel of \cite{christ}, we conclude that
there exist constants $C,D>0$ such that
$\left |\Pi_k\left (\mu_{g}\circ \gamma _X^{-1}(y),y\right)\right|<C\,e^{-D\,k^{1/10}}$, for all $g\in G$,
$k\in \mathbb{N}$
and $y\in \pi^{-1}\left (S_k'\right)$.
The statement follows in view of (\ref{eqn:equivariant-projection}).
\hfill Q.E.D.

\medskip

Given Lemma \ref{lem:annullamento-primo-termine-s},
\begin{eqnarray}\label{eqn:eqz-per-traccia-S_k}
%\lefteqn{}
\mathrm{trace}\left(\Psi_{\varpi,k}\right)\sim\int _{\pi^{-1}(S_k)}\sigma_k(y)\,\Pi_{\varpi,k}\left (\gamma_X^{-1}(y),y\right)\,f(y)\,\mathrm{dens}_{X}(y).
\end{eqnarray}

We next concentrate the integral on progressively shrinking neighborhoods of
$\Phi^{-1}(0)$.
More precisely, we set
\begin{eqnarray}\label{eqn:def-di-t-k}
T_k&=:&\left\{m\in S_k:\mathrm{dist}_M\left (m,\Phi^{-1}(0)\right)<2\,k^{-1/3}\right\},\\
T'_k&=:&\left\{m\in S_k:\mathrm{dist}_M\left (m,\Phi^{-1}(0)\right)>k^{-1/3}\right\},\nonumber
\end{eqnarray}
and let $\{\tau_k,\tau_k'\}$ be a partition of unity on $S_k$ subordinate to the open cover
$\{T_k,T_k'\}$. We may assume that
$\tau_k(m)=\tau\left (\sqrt[3]{k}\,\mathrm{dist}_M\left (m,\Phi^{-1}(0)\right)\right)$,
for a fixed smooth function $\tau:\mathbb{R}\rightarrow \mathbb{R}$.

We shall write $\tau_k$ for $\tau_k\circ \pi$ and similarly for
$\tau_k'$.
Again, insertion of the equality $\tau_k+\tau_k'=1$ in (\ref{eqn:eqz-per-traccia-S_k}) splits
the integral in the sum of two terms, one of which is rapidly decaying as $k\rightarrow +\infty$:

\begin{prop}
\label{prop:annullamento-primo-termine-t}
$\Pi_{\varpi,k}\left (\gamma_X^{-1}(y),y\right)=O\left(k^{-\infty}\right)$
uniformly for $y\in \pi^{-1}(T'_k)$.
\end{prop}

Before commencing the proof of Proposition \ref{prop:annullamento-primo-termine-t},
let us notice that it implies:

\begin{cor}
\label{cor:annullamento-primo-termine-t}
As $k\rightarrow +\infty$, we have
$$
\int _{\pi^{-1}(T'_k)}\big(\sigma_k\cdot \tau_k'\big)(y)\,\Pi_{\varpi,k}\left (\gamma_X^{-1}(y),y\right)\,f(y)\,\mathrm{dens}_{X}(y)
=O\left (k^{-\infty}\right).
$$
\end{cor}

\textit{Proof of Proposition \ref{prop:annullamento-primo-termine-t}.}
If $m\in S_k$ and $k\gg 0$,
there exists a unique $g(m)\in G$, smoothly depending on $m$,
such that
$\mathrm{dist}_M\left(\mu_{g(m)}\circ \gamma^{-1}(m),m\right)=
\mathrm{dist}_M\left(G\cdot \gamma^{-1}(m),G\cdot m\right)$.
If $y\in \pi^{-1}(S_k)$,
 we shall write $g(y)$ for $g\big (\pi(y)\big)$.

We can operate the change of variables $g\rightsquigarrow g\,g(y)$ in
(\ref{eqn:equivariant-projection}), and obtain
\begin{equation}
\Pi_{\varpi,k}\left (\gamma_X^{-1}(y),y\right)=\dim(V_\varpi)\,
\int _G\chi _\varpi\big (g g(y)\big)\,\Pi_k\left (\mu_{g g(y)}\circ \gamma _X^{-1} (y),y\right)\,
d\nu(g).
\label{eqn:equivariant-projection-translated}
\end{equation}

\begin{lem}
\label{lem:aggiusta-le-costanti}
There exist $C,D>0$ such that the following holds: If $m\in S_k$ and
$\mathrm{dist}_G(g,e)>Ck^{-2/5}$, then
$$
\mathrm{dist}_M\left (\mu_{g g(m)}\circ \gamma^{-1}(m),m\right)>D\,k^{-2/5}.
$$
\end{lem}

\noindent
Here $\mathrm{dist}_G$ is the Riemannian distance function on $G$, and $e\in G$ is the unit.

\medskip

\textit{Proof of Lemma \ref{lem:aggiusta-le-costanti}.}
For $m\in \overline{V}$, the map $g\in G\mapsto \mu_g(m)\in M$ is an embedding; by compactness,
there exists $C'>0$ such that
$\mathrm{dist}_M\left (\mu_g(m),m\right)>C'\,\mathrm{dist}_G\left (g,e\right)$ for all $m\in \overline{V}$
and $g\in G$.

Suppose then $m\in S_k$ and $\mathrm{dist}_G(g,e)>Ck^{-2/5}$ for a certain $C>0$. By definition of
$S_k$ and $g(m)$, we have
$\mathrm{dist}_M\left (\mu_{g(m)}(m),m\right)<2\,k^{-2/5}$.
Thus, by the triangle
inequality,
\begin{eqnarray}
\label{eqn:triangle-inequality}
\lefteqn{\mathrm{dist}_M\Big (\mu_{g g(m)}\circ \gamma^{-1}(m),m\Big)}\nonumber\\
&\ge&\mathrm{dist}_M\Big (\mu_{g g(m)}\circ \gamma^{-1}(m),\mu_{g(m)}\circ \gamma^{-1}(m)\Big)-
\mathrm{dist}_M\Big (\mu_{g(m)}\circ \gamma^{-1}(m),m\Big)\nonumber \\
&>& C'\,\mathrm{dist}_G\left (g,e\right)-2\,k^{-2/5}\ge \big (C'\,C-2\big)\,k^{-2/5}.
\end{eqnarray}

\noindent
Given (\ref{eqn:triangle-inequality}), we need only choose $C>2/C'$, $D=C\,C'-2$.
\hfill Q.E.D.

\medskip

As in the proof of Lemma \ref{lem:annullamento-primo-termine-s},
we may then apply the off-diagonal estimates of \cite{christ} to conclude:

\begin{cor}
\label{cor:integra-vicino-e}
Let $C>0$ be as in the statement of Lemma \ref{lem:aggiusta-le-costanti}.
Then $$\Pi_{k}\Big (\mu_{g g(y)}\circ \gamma_X^{-1}(y),y\Big)=O\left (k^{-\infty}\right),$$
uniformly for $y\in \pi^{-1}(S_k)$ and $g\in G$ satisfying $\mathrm{dist}_G(g,e)>Ck^{-2/5}$.
\end{cor}

Let us now set
$$
G_k=:\left\{g\in G:\mathrm{dist}_G\left (g,e\right)<2\,C\,k^{-2/5}\right\},\,
G_k'=:\left\{g\in G:\mathrm{dist}_G\left (g,e\right)>C\,k^{-2/5}\right\},$$
and let
$\{\gamma_k,\gamma_k'\}$ be a smooth partition of unity on $G$ subordinate to the open
cover $\{G_k,G_k'\}$.

Let $\exp_G:\mathfrak{g}\rightarrow G$ be the exponential map, and let
$E\subseteq \mathfrak{g}$ be an open neighborhood of $0$ which is mapped diffeomorphically
under $\exp_G$ to an open neighborhood $U=\exp_G(E)$ of $e$.
Since $G_k\Subset U$ for $k\gg 0$, we may view $\gamma_k$ as a real valued smooth map
on $\mathfrak{g}$ supported on $E$.
With this interpretation, we may assume that $\gamma _k(\xi)=\gamma\left (k^{2/5}\,\xi\right)$
($\xi\in \mathfrak{g}$), for a certain fixed smooth function $\gamma$ on $\mathfrak{g}$.

Inserting the relation $\gamma_k+\gamma_k'=1$ in (\ref{eqn:equivariant-projection}),
integration over $G$ splits as the sum of two terms. In the summand containing $\gamma_k'$,
integration is over $G_k'$; therefore, by Corollary \ref{cor:integra-vicino-e},
if $y\in \pi^{-1}(S_k)$
the integrand is uniformly $O\left(k^{-\infty}\right)$.
Hence we need only worry about the summand containing $\gamma_k$.

To prove Proposition \ref{prop:annullamento-primo-termine-t}, we are thus reduced to proving that
uniformly on $y\in \pi^{-1}(T'_k)$ we have
\begin{equation}
\label{eqn:integrale-vicino-e-piccolo}
\int _{G_k}\gamma_k(g)\,\chi_\varpi\big (g g(y)\big)\,\Pi_k\left (\mu_{g g(y)}\circ \gamma_X^{-1}(y),y\right)\,d\nu(g)=
O\left (k^{-\infty}\right).
\end{equation}
To this end, we shall now invoke the parametrix of the Szeg\"{o} kernel produced in \cite{bs},
and apply an integration by parts as in the proof of the stationary phase Lemma.

Let us recall that the Szeg\"{o} kernel on $X\times X\subseteq   L^*\times L^*$
may be microlocally represented as a Fourier integral operator
of the form
\begin{equation}
\label{eqn:szego-parametrix-bs}
\Pi(x,y)=\int _0^{+\infty}e^{it\psi(x,y)}\,s(x,y,t)\,dt,
\end{equation}
where the complex phase has prescribed Taylor expansion along the diagonal of $L^*\times L^*$,
and the amplitude is a semiclassical symbol admitting an asymptotic expansion
$s(x,y,t)\sim \sum_{j\ge 0}t^{\mathrm{d}-j}\,s_j(x,y)$ \cite{bs}.

Following \cite{z} and \cite{sz},
let us take Fourier components, and perform the change of variable
$t\rightsquigarrow kt$, so that
the left hand side of (\ref{eqn:integrale-vicino-e-piccolo})
may be rewritten
\begin{eqnarray}
\label{eqn:integrale-vicino-e-piccolo-bs}
\lefteqn{\frac{k}{2\pi}\int _{G_k}\int _0^{+\infty}\int _{-\pi}^\pi\gamma_k(g)\,\chi_\varpi\left (g g(y)\right)}\\
&&\cdot e^{
ik\big[t\psi\left (\mu_{g g(y)}\circ r
_{e^{i\theta}}\circ \gamma_X^{-1}(y),y\right)-\theta\big]}\,s\left (\mu_{g g(y)}\circ r
_{e^{i\theta}}\circ \gamma_X^{-1}(y),y,kt\right)\,d\nu(g)\,dt
\,d\theta\nonumber\\
\lefteqn{=\frac{k}{2\pi}\int _{\mathfrak{g}}\int _0^{+\infty}\int _{-\pi}^\pi\gamma\left(k^{2/5}\,\xi\right)\,\chi_\varpi\left (e^{\xi}\,g(y)\right)}\nonumber\\
&&\cdot e^{
i\,k\big[t\psi\left (\mu_{e^{\xi}g(y)}\circ r
_{e^{i\theta}}\circ \gamma_X^{-1}(y),y\right)-\theta\big]}\,s\left (\mu_{e^{\xi}g(y)}\circ r
_{e^{i\theta}}\circ \gamma_X^{-1}(y),y,kt\right)\,H_G(\xi)\,d\xi\,dt\,d\theta \noindent .
\nonumber
\end{eqnarray}
In the latter expression, integration over $G_k$ has been written as an integral over the Lie algebra
$\mathfrak{g}$ by the exponential map $\exp_G(\xi)=e^\xi$, and $H_G(\xi)\,d\xi$ is the pull-back to $\mathfrak{g}$
of the Haar measure on $G$ by $\exp_G$.
Now
(\ref{eqn:integrale-vicino-e-piccolo-bs}) is an oscillatory integral, with phase
$$\Psi(\xi,t,\theta,y)=:t\psi\left (\mu_{e^{\xi}g(y)}\circ r
_{e^{i\theta}}\circ \gamma_X^{-1}(y),y\right)-\theta,$$
depending parametrically on $y$.

If $\xi=0$ and $\gamma \big(\pi(y)\big)=\pi(y)$, so that $g(y)=e$, then
$\Psi=it\left (1-e^{i(\theta+\theta_0)}\right)-(\theta+\theta_0)$. The latter phase was considered
in \cite{z}, \cite{sz};
in this case, $\left |\frac{\partial \Psi}{\partial \theta}\right|>\frac 12$ when $t<\frac 12$.

If more generally $y\in \pi^{-1}(S_k)$ and $g\in G_k$, then $$\mathrm{dist}_M\left (\mu_{gg(y)}\circ r
_{e^{i\theta}}\circ \gamma_X^{-1}(y),y\right)\lesssim k^{-2/5};$$ therefore, by continuity for $k\gg 0$ and $t<\frac 12$
we have
$\left |\frac{\partial \Psi}{\partial \theta}\right|>\frac 13$, say. Hence the contribution from the
locus $t<\frac 12$ is $O\left (k^{-\infty}\right)$.

Similarly, in view of the arguments in \S 3 of \cite{sz},
one can see that
in the same range the contribution coming from
$t\ge 4$, say, is rapidly decreasing.

On the upshot, after disregarding a rapidly decaying contribution, we are left with the oscillatory
integral:
\begin{eqnarray}
\label{eqn:integrale-vicino-e-piccolo-bs-1-2}
\lefteqn{\frac{k}{2\pi}\int _{\mathfrak{g}}\int _{1/2}^{4}\int _{-\pi}^\pi e^{
i k\Psi(\xi,t,\theta,y)}
}\\
&&\cdot \gamma\left(k^{2/5}\,\xi\right)\,\chi_\varpi\left (e^{\xi}\,g(y)\right)\,s\left (\mu_{e^{\xi}g(y)}\circ r
_{e^{i\theta}}\circ \gamma_X^{-1}(y),y,kt\right)\,H_G(\xi)\,d\xi\,dt\,d\theta \nonumber\\
&=&\frac{k}{2\pi}\int _{\mathfrak{g}}\int _{1/2}^{4}\int _{-\pi}^\pi e^{
i k\Psi(\xi,t,\theta,y)}\,\gamma\left(k^{2/5}\,\xi\right)\,S(\xi,kt,\theta,y)\,d\xi\,dt\,d\theta ;\nonumber
\end{eqnarray}
where $S$ is obviously defined;
integration in $\xi$ is supported on a ball centered at $0\in \mathfrak{g}$ and of radius $\sim k^{-2/5}$.

Let us now focus on the directional derivative of $\Psi$ with respect to
$\xi\in \mathfrak{g}$.
This is $\partial _{\xi}\Psi=t\,\partial_{\xi_X} \psi$,
where
$\xi_X\in \mathfrak{X}(X)$ is the vector field generated by $\xi$.
Recall from \cite{bs} that for any $x\in X$,
the differential of $\psi\in \mathcal{C}^\infty(X\times X)$ at
$(x,x)$ is $d_{(x,x)}\psi=(\alpha_x,-\alpha_x)$; more generally,
for any $x\in X$ and $e^{i\theta_0}\in S^1$ we have
\begin{equation}
\label{eqn:differenziale-luno-diag-satura}
d_{ (e^{i\theta_0}x,x)}\psi=\left (e^{i\theta_0}\alpha_{e^{i\theta_0}x},-e^{i\theta_0}\,\alpha_x\right).
\end{equation}

Now if $y\in \pi^{-1}(S_k)$ there exists a unique $e^{i\theta(y)}\in S^1$
such that
$$\mathrm{dist}_X\Big (\mu_{g(y)}\circ r_{e^{i\theta(y)}}\circ \gamma _X^{-1}(y),y\Big)
=\mathrm{dist}_M\Big (\mu_{g(y)}\circ \gamma^{-1}\circ\pi(y),
\pi(y)\Big)\le 2\,k^{-2/5}.
$$
Therefore, there exists $D'>0$ such that for all $y\in \pi^{-1}(S_k)$ and $g\in G_k$
\begin{equation}
\label{eqn:estimate-distance-from-the-diagonal}
\mathrm{dist}_X\Big (\mu_{g g(y)}\circ r_{e^{i\theta(y)}}\circ \gamma_X^{-1}(y),y\Big)\le D'\,k^{-2/5}.
\end{equation}
It follows from (\ref{eqn:differenziale-luno-diag-satura}) and (\ref{eqn:estimate-distance-from-the-diagonal})
that if $(x,y)\in S_k$ and $g\in G_k$, then
\begin{equation}
\label{eqn:diff-di-psi-su-S-k-G-k}
d_{(\mu_{g g(y,x)}\circ \gamma_X^{-1}(y),y)}\psi=\left (e^{-i\theta (y)}\,\alpha
_{\mu_{g g(y)}\circ \gamma_X^{-1}(y)},-e^{-i\theta (y)}\,\alpha_y\right)+O\left(k^{-2/5}\right).
\end{equation}

Let us now make use of the assumption $y\in \pi^{-1}(T_k')$. Since $0\in \frak{g}^*$ is a regular value of $\Phi$,
perhaps after restricting $V$ there exists a constant $E'>0$ such that
$\|\Phi(m)\|\ge E'\,\mathrm{dist}_M\left (m,\Phi^{-1}(0)\right)$, $\forall\, m\in V$.
Therefore, since $\mathrm{dist}_M\left (\cdot,\Phi^{-1}(0)\right)$ is invariant under the $G$-action and $\gamma$,
for $y\in \pi^{-1}(T_k')$ we have
$$\Big \|\Phi\Big(\mu_{g(y)}\circ\gamma ^{-1}\circ\pi(y)\Big)\Big\|\ge E'\,k^{-1/3}.$$
Equivalently, $\forall \,y\in \pi^{-1}(T_k')$
there exists $\eta =\eta(y)\in \mathfrak{g}$ of unit length such that
$\Phi^\eta=:\left <\Phi,\eta\right>$ satisfies
$\left |\Phi^\eta\big(\mu_{ g(y)}\circ\gamma ^{-1}\circ\pi(y)\big)\right|\ge E'\,k^{-1/3}$.
Recalling the definition of $G_k$, setting $E=:E'/2$, say, and letting
$k\gg 0$, we have:

\begin{lem} \label{lem:shrinking-bound}
There exists $E>0$ such that for $k\gg 0$ the following holds:
$\forall \,y\in \pi^{-1}(T_k')$, there exists $\eta =\eta(y)\in \mathfrak{g}$ of unit length such that
$\Phi^\eta=:\left <\Phi,\eta\right>$ satisfies
\begin{equation*}
%\label{eqn:shrinking-bound}
\left |\Phi^\eta\Big(\mu_{g g(y)}\circ \gamma^{-1}\circ\pi(y)\Big)\right|\ge E\,k^{-1/3}.
\end{equation*}
for all $g\in G_k$.
\end{lem}

Recall that
$
\Phi^\eta=-\alpha (\eta_X),
$
where $\eta_X\in \mathfrak{X}(X)$ denotes the smooth vector field generated
by $\eta\in \mathfrak{g}$. Given this and (\ref{eqn:diff-di-psi-su-S-k-G-k}), we conclude
the following: If $ y\in \pi^{-1}(T'_k),\,e^\xi\in G_k$, and $\eta\in \mathfrak{g}$ is as in Lemma \ref{lem:shrinking-bound},
then
\begin{eqnarray}\label{eqn:bound-on-dpsi}
 \left |\partial _{\eta}\Psi\right|&=&t\,\left |\partial _{\eta_X}\psi\right|=t\left |\alpha _{\mu_{e^\xi g(y)}\circ
 \gamma_X^{-1}(y)}\big(\eta_X\big)\right| +O\left (k^{-2/5}\right)\\
 &=&
    t\Big |\Phi^\eta\Big(\mu_{e^\xi g(y )}\circ
 \gamma^{-1}\circ\pi(y)\Big)\Big|+O\left (k^{-2/5}\right)\ge \frac{E}{2}\,k^{-1/3}+O\left (k^{-2/5}\right)\ge \frac E3 \,k^{-1/3},\nonumber
\end{eqnarray}
for all $k\gg 0$, since we are assuming $t\ge \frac 12$.

Let $\{\eta_j\}$ be an orthonormal basis of $\mathfrak{g}$. By
(\ref{eqn:bound-on-dpsi}), for every $ y\in \pi^{-1}(T'_k)$ there exists $j$ such that
$\left |\frac{\partial \Psi}{\partial \eta_j}(\xi,t,\theta,y)\right|>\frac{E}{2\mathrm{g}}\,k^{-1/3}$, whenever
$e^\xi\in G_k$.
In other words, if for every $j=1,\ldots,\mathrm{g}$ we set $s_j=:\frac{\partial \Psi}{\partial \eta_j}$ and
\begin{eqnarray}
\label{eqn:defn-di-v-j}
\lefteqn{V_j=:\Big\{ y\in \pi^{-1}(T'_k):\,}\\ &&\big|s_j(\xi,t,\theta,y)\big|>\frac{E}{2\mathrm{g}}\,k^{-1/3},\,
\forall\,\xi\in \exp_G^{-1}(G_k),t\in \left[1/2,4\right],\,\theta\in [0,2\pi]\Big\},\nonumber
\end{eqnarray}
then $\{V_j\}$ is an open cover of $\pi^{-1}(T'_k)$. Let $\{\varrho_j\}$
be a partition of unity subordinate to this cover; the differential operator on $\pi^{-1}(T'_k)\times \mathfrak{g}$
$$
L=:\sum_{j=1}^\mathrm{g}\left(\frac{\varrho_j}{s_j}\right)\,\frac{\partial}{\partial \eta_j}
$$
satisfies $L(\Psi)=1$, hence $L\left(e^{ik\Psi}\right)=i\,k\,e^{ik\Psi}$.
Recall that in (\ref{eqn:integrale-vicino-e-piccolo-bs-1-2})
integration in $d\xi$ is compactly supported; let us iteratively integrate by parts , so as to obtain
\begin{eqnarray}
\label{eqn:integration-by-parts}
\lefteqn{\int_\mathfrak{g}
e^{
i k\Psi}\,\gamma\left(k^{2/5}\,\xi\right)\,S\,d\xi=\left(\frac{-i}{k}\right)\,\int_\mathfrak{g}
L\left(e^{
i k\Psi}\right)\,\gamma\left(k^{2/5}\,\xi\right)\,S\,d\xi}\nonumber\\
&=&\left(\frac{i}{k}\right)\,\sum_{j=1}^\mathrm{g}\,\int_\mathfrak{g}e^{
i k\Psi}\frac{\partial}{\partial \eta_j}\left(\frac{\varrho_j}{s_j}\,\gamma\left(k^{2/5}\,\xi\right)\,S\right)\,d\xi
\nonumber\\
&=&\left(\frac{i}{k}\right)^2\,\sum_{j_1,j_2=1}^\mathrm{g}\,\int_\mathfrak{g}e^{
i k\Psi}\frac{\partial}{\partial \eta_{j_2}}\left(\frac{\varrho_{j_2}}{s_{j_2}}\,\frac{\partial}{\partial \eta_{j_1}}\left(\frac{\varrho_{j_1}}{s_{j_1}}\,\gamma\left(k^{2/5}\,\xi\right)\,S\right)\right)\,d\xi
\nonumber\\
&=&\cdots \nonumber\\
&=&\left(\frac{i}{k}\right)^r\,\sum_{j_1,\cdots,j_r=1}^\mathrm{g}\,\int_\mathfrak{g}e^{
i k\Psi}Y_J\left(\gamma\left(k^{2/5}\,\xi\right)\,S\right)\,d\xi
\end{eqnarray}
where for any multiindex $J=(j_1,\cdots,j_r)$ and any smooth function $\upsilon$ we have set
$$
Y_J(\upsilon)=:\frac{\partial}{\partial \eta_{j_r}}
\left(
\frac{\varrho_{j_r}}{s_{j_r}}\,
\left(
\frac{\partial}{\partial \eta_{j_{r-1}}}
\left (
\frac{\varrho_{j_{r-1}}}{s_{j_{r-1}}}\,
\left(
\cdots
\frac{\partial}{\partial \eta_{j_1}}
\left(
\frac{\varrho_{1}}{s_{j_1}}\cdot\upsilon
\right)
\cdots
\right)
\right)\right)\right).
$$
For any multindex $B=(b_1,\ldots,b_\mathrm{g})$ let us define
$
s^B=:s_1^{b_1}\cdots s_{\mathrm{g}}^{b_\mathrm{g}}.
$

The following may be proved by induction on $r$:
\begin{lem}
\label{lem:key-estimate-for-decay}
For any $r\in \mathbb{N}$ and $J\in \{1,\ldots,\mathrm{g}\}^r$, we have
\begin{equation}
\label{eqn:develop-for-estimate}
Y_J\Big(\gamma\left(k^{2/5}\,\xi\right)\,S\Big)=\sum _qk^{2a_q/5}\,f_q\,
\cdot\frac{\varrho^{B'_q}}{s^{B_q}},
\end{equation}
where:
\begin{itemize}
  \item $f_q=P_q(s)\cdot\gamma^{(e_q)}\left (k^{2/5}\xi\right)\cdot S^{(f_q)}$, where $P$ is a differential operator
  with no zero order term, and $\gamma^{(e_q)}$, $S^{(f_q)}$ are (possibly) higher order derivatives of
  $\gamma$ and $S$ with respect to $\xi$;   \item $b'_j>0$ if $b_j>0$;
  \item $a_q+|B_q|\le 2r$ for every $q$.
\end{itemize}
\end{lem}

%\textit{Proof.} We proceed by induction on $r$. For $r=1$, the statement is easily verified. Suppose the statement is
%true for $r$, and let us prove it for $r+1$. To fix ideas, let us suppose $j_{r+1}=1$. For each summand appearing in
%(\ref{eqn:develop-for-estimate}), we have

%If $r=1$, we have
%$$
%\frac{\partial}{\partial \eta_j}\left(\frac{\varrho_j}{s_j}\,\gamma\left(k^{2/5}\,\xi\right)\,S\right)=
%\varrho_j\cdot\left (\frac{r_1}{s_j^2}\,S+k^{2/5}\,\frac{r_2}{s_j}\,S+\frac{r_3}{s}\,\frac{\partial S}{\partial
%\eta_j}\right)
%$$
%for bounded functions $r_j$, and the statement is true.

In view of the rescaling $t\rightsquigarrow kt$, the leading term in $S$ and its derivatives grows like $k^{\mathrm{d}}$.
Therefore, the $q$-th summand in Lemma \ref{lem:key-estimate-for-decay}  is bounded
by $$C\,k^{\mathrm{d}+\frac{2}{5}a_q+\frac{1}{3}\,|B_q|}\le C\,k^{\mathrm{d}+\frac{2}{5}(a_q+|B_q|)}
\le C\,k^{\mathrm{d}+\frac{4}{5}r}.$$  In view of (\ref{eqn:integration-by-parts})
we get
\begin{eqnarray}
\label{eqn:key-estimate}
\left|\int_\mathfrak{g}
e^{
i k\Psi}\,\gamma\left(k^{2/5}\,\xi\right)\,S\,d\xi\right|\le C_r\,k^{\mathrm{d}-r/5}
\end{eqnarray}
for any $r\in \mathbb{N}$. This completes the proof of Proposition \ref{prop:annullamento-primo-termine-t}.

\medskip

\hfill Q.E.D.

\medskip

Given Corollary \ref{cor:annullamento-primo-termine-t}, we conclude
\begin{eqnarray}
\label{eqn:final-reduction}
\mathrm{trace}\left (\Psi_{\varpi,k}\right)\sim\int _{\pi^{-1}(T_k)}\varsigma_k(y)\,\Pi_{\varpi,k}\left (\gamma_X^{-1}(y),y\right)\,f(y)
\,\mathrm{dens}_{X}(y),
\end{eqnarray}
where we have set $\varsigma_k=:\sigma_k\cdot \tau_k$.

Let now $\mathrm{Fix}(\gamma_0)\subseteq M_0$ be the fixed locus of $\gamma _0:M_0\rightarrow M_0$,
and set $\widetilde{\mathrm{Fix}(\gamma_0)}=:p^{-1}\left(\mathrm{Fix}(\gamma_0)\right)\subseteq
\Phi^{-1}(0)$, where $p$ is as in (\ref{eqn:defn-di-M-0}).

\begin{prop}
\label{prop:close-to-fixed-locus}
There exists $C>0$ such that
$$
\mathrm{dist}_M\left(G\cdot m,\widetilde{\mathrm{Fix}(\gamma_0)}\right)\le C\,k^{-1/3},\,\,\,\,\,\,\forall\,
m\in T_k.
$$
\end{prop}

\textit{Proof.} If $m\in T_k\subseteq S_k$, then
$\mathrm{dist}_M\big(G\cdot m,G\cdot \gamma(m)\big)<
2\,k^{-2/5}$ by
(\ref{eqn:defn-of-sk}), and
$\mathrm{dist}_M\left (m,\Phi^{-1}(0)\right)<2\,k^{-1/3}$
by (\ref{eqn:def-di-t-k}).

Let $q\in \Phi^{-1}(0)$ be such that
$\mathrm{dist}_M\left (m,\Phi^{-1}(0)\right)=\mathrm{dist}_M\left (m,q\right)$.
Since $\Phi^{-1}(0)$ is $G$-invariant,
$\mathrm{dist}_M\left (m,q\right)=\mathrm{dist}_M(G\cdot m,G\cdot q)$. As $\gamma$ commutes with
the action and preserves the metric, we also have
$$\mathrm{dist}_M\left (m,q\right)=\mathrm{dist}_M\big (\gamma(m),\gamma(q)\big)=
\mathrm{dist}_M\big(G\cdot \gamma(m),G\cdot \gamma(q)\big).$$

Since $G$ acts freely and isometrically on $V$,
there is a Riemannian metric on the manifold $V_0=:V/G$ such that the projection
$\widehat{p}:V\rightarrow V_0$ is a Riemannian submersion. Hence,
$\mathrm{dist}_M\left (G\cdot m,G\cdot n\right)=\mathrm{dist}_{V_0}\left(\widehat{p}(m),\widehat{p}(n)\right)$,
$\forall\,m,n\in V$.
By the triangle inequality on $V_0$,
\begin{eqnarray}
\label{eqn:quasi-punti-fissi}
\lefteqn{\mathrm{dist}_M\Big (G\cdot q,G\cdot \gamma(q)\Big)}\nonumber\\
&\le&
\mathrm{dist}_M\Big (G\cdot q,G\cdot m\Big)+\mathrm{dist}_M\Big (G\cdot m,G\cdot \gamma(m)\Big)
+\mathrm{dist}_M\Big (G\cdot \gamma(m),G\cdot \gamma(q)\Big)\nonumber\\
&=&2\,\mathrm{dist}_M\Big (m,G\cdot \Phi^{-1}(0)\Big)+\mathrm{dist}_M\Big (G\cdot m,G\cdot \gamma(m)\Big)\nonumber\\
&\le&4\,k^{-1/3}+2\,k^{-2/5} < 5\,k^{-1/3}
\end{eqnarray}
if $k\gg 0$.

Set $q_0=:p(q)\in M_0$. Since
$\mathrm{dist}_M\left (q,\widetilde{\mathrm{Fix}(\gamma_0)}\right)=
\mathrm{dist}_{M_0}\left(q_0,\mathrm{Fix}(\gamma_0)\right)$,
\begin{eqnarray}
\label{eqn:inserisci-m'}
\mathrm{dist}_M\left (m,\widetilde{\mathrm{Fix}(\gamma_0)}\right)
&\le&\mathrm{dist}_M\left (m,q\right)+\mathrm{dist}_M\left (q,\widetilde{\mathrm{Fix}(\gamma_0)}\right)\nonumber\\
&\le&2\,k^{-1/3}+\mathrm{dist}_{M_0}\left(q_0,\mathrm{Fix}(\gamma_0)\right).
\end{eqnarray}

\begin{lem}
\label{lem:q-0-vicino-luogo-fisso}
There exists a constant $C>0$ such that for all $k\gg 0$ we have
 $$\mathrm{dist}_{M_0}\left(q_0,\mathrm{Fix}(\gamma_0)\right)\le
C\,k^{-1/3}.$$
\end{lem}

\textit{Proof.}
By (\ref{eqn:quasi-punti-fissi}), for all $k\gg 0$ we have
\begin{eqnarray}
\label{eqn:vicino-luogo-fisso-1}
\mathrm{dist}_{M_0}\Big (q_0,\gamma_0(q_0)\Big)=\mathrm{dist}_M\Big (G\cdot q,G\cdot \gamma(q)\Big)<
5\,k^{-1/3}.
\end{eqnarray}

Let now $F_1,\ldots, F_\ell\subseteq M_0$, with normal bundles $N_1,\ldots,N_l$, be as in Definition
\ref{defn:componenti-luogo-fisso}.
Let $\exp_{l}:N_l\rightarrow M_0$ be the exponential map, $(q_0',n)\mapsto \exp_{l}(q_0',n)=:\exp_{q_0'}(n)$.
For $\epsilon>0$, let $N_{l}^{(\epsilon)}=:\left \{(q_0',n)\in N_l:\|n\|<\epsilon\right\}$.
Choose $\epsilon>0$ so small that $\exp_l$ induces a diffeomorphism between
$N_{l}^{(\epsilon)}$ and an open neighborhood $F_{l}^{(\epsilon)}\subseteq M_0$ of $F_l$, and
$\overline{F_{l_1}^{(\epsilon)}}\cap \overline{F_{l_2}^{(\epsilon)}}=\emptyset,\,\forall\,l_1\neq l_2\in \{1,\ldots,\ell\}$.

If $q_0'\in F_l$, the normal exponential map
$\exp_{q_0'}:N_{l,q_0'}\rightarrow M_0$ is an isometric immersion
at the origin. By compactness of $F_l$,
perhaps after decreasing $\epsilon$ we may assume that if
$q_0'\in F_l$ and $(q_0',n),(q_0',n')\in N_{l,q_0'}\cap N_{l,\epsilon}$ then
\begin{equation}
\label{eqn:exp-almost-isom}
2\,\|n-n'\|\ge \mathrm{dist}_{M_0'}\Big (\exp_{N_l}(q_0',n),\exp_{N_l}(q_0',n')\Big)\ge \frac 12\,\|n-n'\|.
\end{equation}

There exists $\delta>0$ such that
$$
\mathrm{dist}_{M_0}\Big (q_0,\mathrm{Fix}(\gamma_0)\Big)\ge \epsilon\,
\Rightarrow\,\mathrm{dist}_{M_0}\Big (q_0,\gamma_0(q_0)\Big)\ge \delta.
$$
Thus, if $k\gg 0$ and
(\ref{eqn:vicino-luogo-fisso-1}) holds, then $q_0\in \bigcup _{l=1}^\ell F_{l}^{(\epsilon)}$; hence,
$q_0=\exp_{N_l}(q_0',n)$ for some $(q_0',n)\in N_{l}^{(\epsilon)}$.
Given that $\gamma_0:M_0\rightarrow M_0$ is a Riemannian isometry, we have
$\exp _{M_0}\circ d\gamma_0=\gamma _0\circ \exp_{M_0}:TM_0\rightarrow M_0$.
In view of (\ref{eqn:exp-almost-isom}), we deduce
\begin{eqnarray}
\label{eqn:stime-dianza-autovalori}
\lefteqn{\mathrm{dist}_{M_0}\Big (q_0,\gamma_0(q_0)\Big)=
\mathrm{dist}_{M_0}\Big (\exp_{M_0}(q_0',n),\gamma_0\circ \exp_{M_0}(q_0',n)\Big)}\\
&=& \mathrm{dist}_{M_0}\Big (\exp_{M_0}(q_0',n),\exp_{M_0}\circ d_{q_0'}\gamma_0(n)\Big)\ge \frac 12\,\left\|
 d_{q_0'}\gamma_0(n)-n\right\|\nonumber\\
&\ge&  \frac 12\,\inf \big\{\left |\lambda_i-1\right|\big\}\,\|n\| \ge \frac 14\,
\inf \big\{\left |\lambda _i-1\right|\big\}\,\mathrm{dist}_{M_0}\Big (q_0,\mathrm{Fix}(\gamma_0)
\Big).\nonumber
\end{eqnarray}
Since $\lambda_i\neq 1$ for all $i$,
the statement follows from (\ref{eqn:vicino-luogo-fisso-1}) and (\ref{eqn:stime-dianza-autovalori}).
\hfill Q.E.D.

\medskip

Proposition \ref{prop:close-to-fixed-locus} now follows from (\ref{eqn:inserisci-m'})
and Lemma \ref{lem:q-0-vicino-luogo-fisso}.

\hfill Q.E.D.

\medskip

It is now in order to give the:

\textit{Proof of Lemma \ref{lem:inotorno-phi-0}.} This follows from a simplified version of the previous
arguments. Recall that
\begin{equation}
\Pi_{\varpi,k}\left (x,y\right)=
\dim(V_\varpi)\,\int _G\chi _\varpi\left (g\right)\,\Pi_k\big (\mu_g(x),y\big)\,
d\nu(g),
\label{eqn:equivariant-projection-general}
\end{equation}
for all $(x,y)\in X\times X$.

If $R\subseteq X\times X\setminus I(\Phi)$ is compact,
there exists $\delta>0$ such that $\forall\,
(x,y)\in R$ we have
$$
\max\Big\{\mathrm{dist}_M\big(G\cdot \pi(x),G\cdot \pi(y)\big),\mathrm{dist}_M\left(\pi(x),\Phi^{-1}(0)\right)\Big\}
\ge \delta.
$$
The case $\mathrm{dist}_M\big(G\cdot \pi(x),G\cdot \pi(y)\big)\ge \delta$ can be
handled by the arguments in the proof
of Lemma \ref{lem:annullamento-primo-termine-s}, replacing the lower bound $k^{-2/5}$ used there with
$\delta$.

To deal with the case $\mathrm{dist}_M\left(\pi(x),\Phi^{-1}(0)\right)>\delta$, we may by the same argument
restrict to the case where $\mathrm{dist}_M\big(G\cdot \pi(x),G\cdot \pi(y)\big)<\epsilon$, for some
fixed but arbitrarily small $\epsilon >0$. Invoking the fact that the Szeg\"{o} kernel is smoothing away from
the diagonal, we conclude that we only miss a rapidly decaying contribution if we restrict the
$G$-integration in (\ref{eqn:equivariant-projection-general}) to the open subset
$G(x,y)=:\left\{g\in G:\mathrm{dist}_M\big(\mu_g\circ \pi(x),\pi(y)\big)<2\epsilon \right\}$, say
(the introduction of an appropriate partition of unity on $G$ is understood). Now the hypothesis implies that
$\big\|\Phi\circ \pi(x)\big\|>\delta'$ for some $\delta'\gtrsim \delta$. By the arguments leading to
the proof of Lemma \ref{lem:shrinking-bound} and (\ref{eqn:bound-on-dpsi}), it follows
that if $0<\epsilon\ll 1$ then $\big|\partial _\eta\Psi\big|>\delta'/4$ on the range of integration.
The statement then follows by a simpler version of the argument following (\ref{eqn:bound-on-dpsi}).

\hfill Q.E.D.

\medskip

Summing up, $T_k$ is a shrinking open neighborhood of
$\widetilde{\mathrm{Fix}(\gamma_0)}$.
To obtain an asymptotic expansion for $\mathrm{trace}(\Psi_{\varpi,k})$,
we shall insert in (\ref{eqn:final-reduction}) the scaling limit
asymptotics for $\Pi_{\varpi,k}$ proved in \cite{p}.
Scaling asymptotics are most naturally stated in local Heisenberg coordinates; therefore, we shall
cover $\pi^{-1}\left(\widetilde{\mathrm{Fix}(\gamma_0)}\right)$
by invariant open sets with a \lq transverse Heisenberg structure\rq,
providing convenient coordinates to perform the integration.

For $r\in \mathbb{N}$ and $\epsilon>0$, let $B_r(\epsilon)$ be the open ball
centered at $\mathbf{0}\in \mathbb{R}^r$ and of radius $\epsilon$.
Referring the reader to \cite{sz} for the precise definitions, we recall that a system of local Heisenberg coordinates
for the circle bundle $X$ centered at a given $x\in X$ is determined by the following data: i) a preferred local chart
$\mathfrak{f}:B_{2\mathrm{d}}(\epsilon)\rightarrow U\subseteq M$ for the
(almost) K\"{a}hler manifold $(M,\omega,J)$ centered at $m=\pi(x)$, and
ii) a preferred local frame $e:U\rightarrow A$ at $m$ satisfying $e^*(m)=x$, where $\langle e^*,e\rangle =1$.
That $\mathfrak{f}$ is a preferred local chart at $m$ means that it trivializes the unitary structure of $T_mM$.
In the present integrable setting, $\mathfrak{f}$ may be chosen holomorphic, but this
is not necessary, and won't be assumed in the following; $\mathfrak{f}$ is at any rate always holomorphic and symplectic at
$\mathbf{0}\in B_{2\mathrm{d}}(\epsilon)\subseteq \mathbb{C}^\mathrm{d}$.
In fact,
upon choosing an orthonormal complex basis of $T_mM$, the exponential map
$\exp_m:\mathbb{C}^\mathrm{d}\cong T_mM\rightarrow M$ restricts to a preferred local
chart on $B_{2\mathrm{d}}(\epsilon)$, for some $\epsilon>0$; at places
it will simplify our arguments to make this choice.
Explicitly, given $\mathfrak{f}$ and $e$ the associated Heisenberg local chart is then
$$
\psi:B_{2\mathrm{d}}(\epsilon)\times (-\pi,\pi)\rightarrow
\pi^{-1}(U),\,\,\,\,\,\,(\mathbf{z},\vartheta)\mapsto
e^{i\vartheta}\,\frac{e^*\big(\mathfrak{f}(\mathbf{z})\big)}{\left \|e^*\big(\mathfrak{f}(\mathbf{z})\big)\right\|}.
$$
Following \cite{sz},
we set
\begin{equation}
\label{eqn:sz-notation}
x+w=:\psi\big(w,0)
\end{equation}
if $w\in T_mM\cong \mathbb{C}^\mathrm{d}$, $\|w\|<\epsilon$. In this notation, $T_mM$ is implicitly identified
with the horizontal subspace $\mathrm{Hor}_x(X)\subseteq T_xX$ for the connection.

Given any $x\in X$, it is always possible to find local Heisenberg coordinates centered at $x$, and this
construction may be deformed smoothly with $x$; that is, given any $x\in X$ there exist $x\in U\subseteq X$ open and
a smooth map $\Psi:U\times B_{2\mathrm{d}}(\epsilon)\times (-\pi,\pi)\rightarrow X$, such that for any $y\in U$
the partial map $\psi^{(y)}=:\Psi(y,\cdot,\cdot):B_{2\mathrm{d}}(\epsilon)\times (-\pi,\pi)\rightarrow X$ is a Heisenberg
local chart for $X$ centered at $y$. We may, and will, assume without loss that
$$\Psi\left(r_{e^{i\vartheta_0}}(y),\mathbf{z},\vartheta\right)=\Psi\left(y,\mathbf{z},\vartheta
-\vartheta_0\right)$$
whenever the two sides are defined.

In the present equivariant setting, suppose $x\in X$ and let $\psi:B_{2\mathrm{d}}(\epsilon)\times (-\pi,\pi)\rightarrow X$ be a system of local Heisenberg coordinates centered at
$x$, associated to the preferred choices $\mathfrak{f}$ and $e$. Then for any $g\in G$ we obtain a system of Heisenberg coordinates
centered at $\mu_g(x)$ by considering the composition $\psi_g=:\mu_g\circ \psi$. Clearly, $\psi_g$ is associated to the preferred
choices $\mu_g\circ \mathfrak{f}$ and $\widehat{\mu}_g(e)$; here $\widehat{\mu}$ denotes the action on the collection of local sections
of $A$.

Let us set, for ease of notation, $M'=:\Phi^{-1}(0)\subseteq M$, $X'=:\pi^{-1}\left(\Phi^{-1}(0)\right)
\subseteq X$, and let us denote by $\pi':X'\rightarrow M'$ the projection. Then $G$ acts freely on $M'$ and
$X'$, and we have the commutative diagram:
\begin{equation}
\label{eqn:commutative-diagram}
\begin{array}{ccccclc}
   &  &  & \widetilde{p} &  &  &  \\
   &  X'& &\longrightarrow  &  & X_0=X'/G &  \\
   &  &  &  &  &  &  \\
   \pi'& \downarrow &  &  &  & \downarrow & \pi_0 \\
   &  &  &  &  &  &  \\
   & M' &  & \longrightarrow &  & M_0=M'/G &  \\
   &  &  & p &  &  &
\end{array}
\end{equation}
where the vertical arrows are principal $S^1$-bundles and the horizontal arrows are principal
$G$-bundles. $\pi'$ is $G$-equivariant and $\widetilde{p}$ is $S^1$-equivariant.
Given subsets $V\subseteq M_0$ and $U\subseteq M$, we shall set
$X_0(V)=:\pi_0^{-1}(V)\subseteq X_0$, $X(V)=:\pi^{-1}(V)\subseteq X$.

Suppose given:
\begin{itemize}
  \item $m_0\in F_l\subseteq \mathrm{Fix}(\gamma_0)\subseteq M_0$;
  \item an open subset $V\subseteq F_l$ with
$m_0\in V$;
  \item a smooth section $\sigma:V\rightarrow \widetilde{F}_l=:p^{-1}(F_l)$ of the principal $G$-bundle $\widetilde{F_l}
\rightarrow F_l$.
\end{itemize}

\noindent
%Perhaps after restricting $V$, we may assume that $V$ is the image of a coordinate chart
%$\beta:B_{2\mathrm{d}_l}(\epsilon)\rightarrow V$ for $F_l$ centered at $m_0$.
Then there exists a unique smooth section
$\widetilde{\sigma}:X_0(V)\rightarrow X'$ of $\widetilde{p}$ which is a lift of $\sigma$, that is, such that
$\sigma\circ \pi_0=\pi'\circ \widetilde{\sigma}:X_0(V)\rightarrow M'$.
$\widetilde{\sigma}$ is necessarily $S^1$-equivariant.

By the above, we may also suppose given a smooth map
\begin{equation}
\label{eqn:heisenberg-mobile}
\Psi:X_0(V)\times B_{2\mathrm{d}}(\epsilon)\times (-\pi,\pi)\rightarrow X,
\end{equation}
such that for any $x_0'\in X_0(V)$ the partial map
$\psi^{(x_0')}=:\Psi(x_0',\cdot,\cdot)$ is a Heisenberg local chart for $X$ centered at $\widetilde{\sigma}(x_0')$, with image containing
some fixed
open neighborhood $X(U)\supseteq \pi^{-1}\big(\sigma (m_0)\big)$.
We shall write $$\left(\psi^{(x_0')}\right)^{-1}=
\left (z_1^{(x_0')},\ldots,z_\mathrm{d}^{(x_0')},\vartheta^{(x_0')}\right)
:U\rightarrow B_{2\mathrm{d}}(\epsilon)\times (-\pi,\pi)
$$
for the corresponding Heisenberg local chart. Here we identify $\mathbb{R}^{2\mathrm{d}}\cong \mathbb{C}^\mathrm{d}$ in the standard
manner, and
$z_j^{(x_0')}:U\rightarrow \mathbb{C}$ is a smooth function. Let
$a_j^{(x_0')}=:\Re\left(z_j^{(x_0')}\right),b_j^{(x_0')}=:\Im\left(z_j^{(x_0')}\right):
X(U)\rightarrow \mathbb{R}$.
Clearly, $a_j^{(x_0')}$ and $b_j^{(x_0')}$ descend to $U$, and form the system $\mathfrak{f}^{(x_0')}$
of preferred local coordinates on $M$
centered at $\pi\circ \widetilde{\sigma}(x_0')$ which underlies $\psi^{(x_0')}$.

After composing with a suitable local diffeomorphism of $M$, smoothly varying with $x_0'$,
we may assume that for every $x_0'\in X_0(V)$ the following conditions are satisfied by
$\mathfrak{f} ^{(x_0')}$:

\begin{enumerate}
  \item $\Phi^{-1}(0)\cap U=\left\{b_{\mathrm{d}-\mathrm{g}+1}^{(x_0')}=\cdots =b_\mathrm{d}^{(x_0')}=0\right\}$;
  \item $\widetilde{F_l}\cap U=\left\{z^{(x_0')}_{\mathrm{d}_l+1}=\cdots=z^{(x_0')}_{\mathrm{d}-\mathrm{g}}=b_{\mathrm{d}-\mathrm{g}+1}^{(x_0')}=\cdots =b_\mathrm{d}^{(x_0')}=0\right\}$.
\end{enumerate}

By the previous discussion, composing with the $G$-action we then obtain a smooth map
(we write $\mu$ for $\mu_X$)
$$
\widetilde{\Psi}:G\times X_0(V)\times B_{2\mathrm{d}}(\epsilon)\times (-\pi,\pi)\rightarrow X,\,\,\,\,\,
(g,x_0',\mathbf{z},\vartheta)\mapsto \mu_g\Big(\Psi\big(x_0',\mathbf{z},\vartheta\big)\Big),
$$
such that for any $(g,x_0')\in G\times X_0(V)$ the partial map
\begin{equation}
\label{eqn:heisenber-local-chart-0}
\psi^{(g,x_0')}=:\widetilde{\Psi}\big(g,x_0',\cdot,\cdot):B_{2\mathrm{d}}(\epsilon)\times (-\pi,\pi)\rightarrow X
\end{equation}
is a Heisenberg local chart for
$X$ centered at $\mu_g\circ \widetilde{\sigma}(x_0')$, and whose image contains $\mu_g\big(X(U)\big)$.
We shall denote by
\begin{eqnarray}
\label{eqn:heisenber-local-chart}
\lefteqn{\left(\psi^{(g,x_0')}\right)^{-1}=\left(\psi^{(x_0')}\right)^{-1}\circ \mu_g^{-1}}\\
&=&\left (z_1^{(g,x_0')},\ldots,z_\mathrm{d}^{(g,x_0')},\vartheta^{(g,x_0')}\right)
:\mu_g(U)\rightarrow B_{2\mathrm{d}}(\epsilon)\times (-\pi,\pi)
\nonumber
\end{eqnarray}
the corresponding Heisenberg local coordinates. By the $G$-invariance of $\Phi^{-1}(0)$ and
$\widetilde{F}_l$ we obtain that for every $(g,x_0')\in G\times V$

\begin{enumerate}
  \item $\Phi^{-1}(0)\cap \mu_g(U)=\left\{b_{\mathrm{d}-\mathrm{g}+1}^{(g,x_0')}=\cdots =b_\mathrm{d}^{(g,x_0')}=0\right\}$;
  \item $\widetilde{F_l}\cap \mu_g(U)=\left\{z^{(g,x_0')}_{\mathrm{d}_l+1}=
  \cdots=z^{(g,x_0')}_{\mathrm{d}-\mathrm{g}}=
  b_{\mathrm{d}-\mathrm{g}+1}^{(g,x_0')}=\cdots =b_\mathrm{d}^{(g,x_0')}=0\right\}$.
\end{enumerate}

Having in mind the identifications
$$
\mathbb{R}^{2\mathrm{d}}
\cong \mathbb{R}^{2\mathrm{d}_l}\times \mathbb{R}^{2\mathrm{c}_l}\times \mathbb{R}^{2\mathrm{g}}\cong
\mathbb{C}^{\mathrm{d}_l}\times \mathbb{C}^{\mathrm{c}_l}\times \mathbb{C}^{\mathrm{g}},
$$
the following is a straightforward consequence of the previous discussion:

\begin{lem}
\label{lem:equivariant-embedding}
Suppose $\epsilon>0$ is sufficiently small, and define
$$\Upsilon:G\times X(V)\times B_{2\mathrm{c}_l}(\epsilon)\times
B_{\mathrm{g}}(\epsilon)\rightarrow X$$
by
$$
\Upsilon\Big(g,x_0',\mathbf{z},\mathbf{b}\Big)=:
\widetilde{\Psi}\Big(g,x_0',\big(\mathbf{0},\mathbf{z},i\mathbf{b},0\big)\Big)=
\psi^{(g,x_0)}\big(\mathbf{0},\mathbf{z},i\mathbf{b},0\big),
$$
where $\mathbf{z}=\big(z_{\mathrm{d}_l+1},\cdots,z_{\mathrm{d}-\mathrm{g}}\big)\in B_{2\mathrm{c}_l}(\epsilon)
\subseteq \mathbb{C}^{\mathrm{c}_l}$, and $\mathbf{0}$ denotes the origin of $\mathbb{C}^{\mathrm{d}_l}$.
Then:
\begin{enumerate}
  \item $\Upsilon$ is an equivariant diffeomorphism onto a $(G\times S^1)$-invariant open neighborhood $B$ of $
  G\cdot \widetilde{\sigma}\big(X_0(V)\big)\subseteq X$.
  \item In terms of $\Upsilon$, $\widetilde{\mathrm{Fix}(\gamma_0)}\cap B$ is defined by the conditions
  $\mathbf{z}=\mathbf{0}$, $\mathbf{b}=\mathbf{0}$; in other words,
  $\Upsilon^{-1}\left (\widetilde{\mathrm{Fix}(\gamma_0)}\right)=G\times X(V)\times \{\mathbf{0}\}\times \{\mathbf{0}\}
  $.
\end{enumerate}
\end{lem}

Before we proceed, let us dwell on the local structure of $M$ along $\widetilde{F}_l=p^{-1}(F_l)\subseteq M'$.
For any $m\in M$, let $\mathfrak{g}_M(m)\subseteq T_mM$ be the vector subspace generated by the infinitesimal action of
$\mathfrak{g}$.
Thus $\mathfrak{g}_M$ is a rank-g vector subbundle of $TM$ on some invariant open neighborhood of $M'$.
If $m\in M'$, we have the unitary direct sum decompositions
\begin{equation}
\label{eqn:struttura-locale-M-phi}
T_mM=J_m\big (\mathfrak{g}_M(m)\big)\oplus T_mM',\,\,\,\,\,\,\,
T_mM'=\mathfrak{g}_M(m)\oplus H_m;
\end{equation}
here $H_m=:T_mM'\cap \big(\mathfrak{g}_M(m)\big)^\perp$ is a complex subspace, that gets unitarily identified
with $T_{p(m)}M_0$ under $d_{m}p$ (the superscript $\perp$ stands for \lq Euclidean orthocomplement\rq).
Thus if $m\in \widetilde{F}_l$ with this identification we also have
\begin{equation}
\label{eqn:struttura-locale-M-phi-hor}
H_m\cong T_{p(m)}M_0=T_{p(m)}F_l\oplus (N_l)_{p(m)}, \,\,\,\,\,\,
T_m\widetilde{F}_l\cong T_{p(m)}F_l\oplus \mathfrak{g}_M(m),
\end{equation}
where $(N_l)_{p(m)}$ denotes the fiber at $p(m)$ of the normal bundle $N_l$ of $F_l\subseteq M_0$.

If $m\in M'$, let us set
\begin{equation}
\label{eqn:struttura-locale-M-phi-vert}(T_mM)_\mathrm{t}=:J_m\big (\mathfrak{g}_M(m)\big),\,\,(T_mM)_\mathrm{v}=:\mathfrak{g}_M(m),\,\,
(T_mM)_\mathrm{h}=:H_m.
\end{equation}
Here the suffix t stands for \lq transverse to $M'$\rq, v stands for \lq vertical for the principal
$G$-bundle structure of $M'\rightarrow M_0$\rq, h for \lq horizontal\rq.

If in addition $m\in \widetilde{F}_l$, with a slight abuse of language, let us set
\begin{equation}
\label{eqn:struttura-locale-M-phi-vert-1}
(T_mM)_\mathrm{h,tg}=:T_{p(m)}F_l,\,\,(T_mM)_\mathrm{h,nor}=:(N_l)_{p(m)}.
\end{equation}
Thus, $(T_mM)_\mathrm{h}=(T_mM)_\mathrm{h,tg}\oplus (T_mM)_\mathrm{h,nor}$.
Here the suffix h,tg stands for \lq horizontal and tangent to $\widetilde{F}_l$\rq,
h,nor for \lq horizontal and normal to $\widetilde{F}_l$\rq.
Accordingly, if $m\in \widetilde{F}_l$ any $v\in T_mM$ may be decomposed as
$v=v_\mathrm{t}+v_\mathrm{v}+v_\mathrm{h,tg}+v_\mathrm{h,nor}$.

Let us consider again the statement of Lemma \ref{lem:equivariant-embedding}.
%Recall that for any $x\in X$ we implicitly identify $T_{\pi(x)}M$ with the horizontal subspace of
%the connection $\mathrm{Hor}_x(X)\subseteq T_xX$.
Suppose $(g,x_0')\in G\times X(V)$.
Since $\widetilde{\Psi}^{(g,x_0')}$ is a local Heisenberg chart for $X$ centered at $x=:\mu_g\circ\widetilde{\sigma}
(g,x_0')$, and satisfying 1. and 2. above, any $\mathbf{z}\in \mathbb{C}^{\mathrm{c}_l}$ gets identified with an
appropriate
$v_\mathrm{h,nor}\in (T_mM)_\mathrm{h,nor}$
where $m=\pi(x)\in \widetilde{F}_l$.
Similarly, if $\mathbf{b}\in \mathbb{R}^\mathrm{g}$ then $i\,\mathbf{b}$ gets identified with an appropriate
$v_\mathrm{t}\in (T_mM)_\mathrm{t}$.
Following (\ref{eqn:sz-notation}), for sufficiently small
$\mathbf{z}\in \mathbb{C}^{\mathrm{c}_l}$ and $\mathbf{b}\in \mathbb{R}^\mathrm{g}$
we then have
$$
\Upsilon\Big(g,x_0',\mathbf{z},\mathbf{b}\Big)=\mu_g\left (\widetilde{\sigma} (x'_0)\right)+v_\mathrm{h,nor}
+v_\mathrm{t}
$$

For
$\epsilon>0$ sufficiently small, we can then replace $B$ in Lemma \ref{lem:equivariant-embedding}
by
\begin{eqnarray}
\label{eqn:b}
\lefteqn{B=:\Big \{\mu_g\left (\widetilde{\sigma} (x'_0)\right)+v_\mathrm{h,nor}
+v_\mathrm{t}:(g,x_0')\in G\times X(V), }\\
&&v_\mathrm{h,nor}\in
\left(T_{\mu_g\circ\pi\circ \widetilde{\sigma} (x'_0)}M\right)_\mathrm{h,nor},
\,v_\mathrm{t}\in \left(T_{\mu_g\circ\pi\circ \widetilde{\sigma} (x'_0)}M\right)_\mathrm{t},
\|v_\mathrm{h,nor}\|,\,\|v_\mathrm{t}\|<\epsilon\Big\}, \nonumber
\end{eqnarray}
an invariant open neighborhood of $\pi^{-1}\left(p^{-1}(V)\right)
\subseteq \pi^{-1}\left(\widetilde{\mathrm{Fix}(\gamma_0)}\right)$.

Let $B_\mathrm{g}(\epsilon)\subseteq \mathbb{R}^\mathrm{g}$,
$B_{2\mathrm{d}_l}(\epsilon)\subseteq \mathbb{C}^{\mathrm{d}_l}$,
$B_{2\mathrm{c}_l}(\epsilon)\subseteq \mathbb{C}^{\mathrm{c}_l}$
the open balls of radius $\epsilon$ centered at the origin.
The paramerization (\ref{eqn:b})
defines a diffeomorphism
\begin{eqnarray}
\label{eqn:parametr-per-tilde-b}
\mathcal{T}: G\times X(V)\times B_\mathrm{g}(\epsilon)\times B_{2\mathrm{c}_l}(\epsilon)&\longrightarrow& B,
\nonumber\\
%\left(e^{i\theta},g,h,x_0',v_\mathrm{t},v_\mathrm{h,nor},w_\mathrm{t},w_\mathrm{h,tg},w_\mathrm{h,nor}\right)&\mapsto&
\mathcal{T}\left(g,x_0',v_\mathrm{t},v_\mathrm{h,nor}\right)&=:&\mu_g\left (\widetilde{\sigma} (x'_0)\right)+v_\mathrm{h,nor}
+v_\mathrm{t}.
\end{eqnarray}

\noindent
%%%%%Furthermore, there is a natural smooth projection
%%%%%$q:\widetilde{B}\rightarrow V\subseteq F_l$, given by
%%%%
%%%%%%%%%\begin{equation}\label{eqn:projection-b-v}
 %%%%%%%   q\Big(\mu_{h}\circ r_{e^{i\theta}}\Big (\mu_g\circ\widetilde{\sigma} (x'_0)
%%%%%%%%%+v\Big),\mu_g\circ\widetilde{\sigma} (x'_0)+v_\mathrm{h,nor}
%%%%%%%+v_\mathrm{t}\Big)=:\pi_0\left (x_0'\right).
%%%\end{equation}
%%%%

Now for every $l=1,\ldots,\ell$ let $\{V_{lj}\}_j$ be a finite open cover of $F_l$,
such that on every $V_{lj}$ there is defined a smooth section $\sigma_{lj}:V_{lj}\rightarrow M'$
of $p$; in particular, $\left\{p^{-1}(V_{lj})\right\}_{l,j}$ is an invariant
open cover of $\widetilde{\mathrm{Fix}(\gamma_0)}$.

For every $l,j$, let $B_{lj}\subseteq X$ defined by (\ref{eqn:b}) with
$\left (V_{lj},\sigma_{lj}\right)$ in place of $(V,\sigma)$; hence
$B_{lj}$ is an $S^1\times G$-invariant open neighborhood of $\pi^{-1}\left(p^{-1}(V_{lj})\right)$.
%Let
%$q_{lj}:\widetilde{B}_{lj}\rightarrow V_{lj}$ be defined accordingly as in (\ref{eqn:projection-b-v}).

Finally, set $E=:\bigcup _{lj}B_{lj}$. Then
$E$ is an $S^1\times G$-invariant open neighborhood of $\pi^{-1}\left(\widetilde{\mathrm{Fix}(\gamma_0)}\right)\subseteq X$,
and $\{B_{lj}\}_{l,j}$ is an
open cover of $E$. Let $\{\tau_{lj}\}_{l,j}$ be a smooth partition of unity on $E$ subordinate to this
cover. After averaging, we may assume that each $\tau_{lj}$ is $S^1\times G$-invariant;
hence each $\tau_{lj}$
descends in a natural manner to a smooth function $\tau_{lj}^0$ on $F_l$, and
$\{\tau_{lj}^0\}_{l,j}$ is a smooth partition of unity on $\mathrm{Fix}(\gamma_0)$, subordinate to
the open cover $\{V_{lj}\}$.

Let us now return to (\ref{eqn:final-reduction}). Since $T_k$
is a shrinking open neighborhood of $\widetilde{\mathrm{Fix}(\gamma_0)}$ as $k\rightarrow +\infty$,
we have $\pi^{-1}(T_k)\subseteq E$ for all
$k\gg 0$, hence
$$\pi^{-1}(T_k)=\bigcup _{lj}\pi^{-1}(T_k)\cap B_{lj}.$$
Thus $\pi^{-1}(T_k)\cap
B_{lj}$ is a shinking open neighborhood of $\pi^{-1}\left (p^{-1}(V_{lj})\right)$.
Inserting the relation $\sum _{l,j}\tau_{lj}=1$ in (\ref{eqn:final-reduction}) yields
$\mathrm{trace}(\Psi_{\varpi,k})\sim \sum _{lj}\mathrm{trace}(\Psi_{\varpi,k})_{lj}$, where
\begin{eqnarray}\label{eqn:trace-lj-th}
%\lefteqn{\mathrm{trace}\left (\Psi_{\varpi,k}\right)_{lj}}
\mathrm{trace}\left (\Psi_{\varpi,k}\right)_{lj}=:\int _{\pi^{-1}(T_k)\cap B_{lj}}\tau_{lj}(y)\,\varsigma_k(y)\,\Pi_{\varpi,k}\left (\gamma_X^{-1}(y),y\right)\,f(y)
\,\mathrm{dens}_{X}(y).
%\nonumber
\end{eqnarray}
Let us now estimate asymptotically each summand
(\ref{eqn:trace-lj-th}). To simplify our notation, in the following
formulae we shall temporarily fix a pair $(l,j)$, and occasionally
write $\sigma$, $\widetilde{\sigma}$ for $\sigma_{lj}$,
$\widetilde{\sigma}_{lj}$.

We can parametrize $\pi^{-1}(T_k)\cap
B_{lj}$ by (\ref{eqn:b}).
Given (\ref{eqn:defn-of-sk}), (\ref{eqn:def-di-t-k}), and Proposition
\ref{prop:close-to-fixed-locus}
$\|v_\mathrm{t}\|,\,\|v_\mathrm{h,nor}\|\lesssim
k^{-1/3}$. Inserting (\ref{eqn:b}) in (\ref{eqn:trace-lj-th}),
we shall rescale $\|v_\mathrm{t}\|,\,\|v_\mathrm{h,nor}\|$ by a factor $k^{-1/2}$ and
integrate the rescaled variables over a ball of radius
$\thickapprox k^{1/6}$. In other words, we shall write
\begin{eqnarray}
\label{eqn:x-e-y-riscalati}
y&=&\mu_g\circ\widetilde{\sigma} (x'_0)+\frac{1}{\sqrt{k}}\,\big(v_\mathrm{t}+v_\mathrm{h,nor}\big),
%%%%y&=&\mu_{h}\circ r_{e^{i\theta}}\left
%%%%%%%%(\mu_g\circ\widetilde{\sigma}(x'_0)
%%%%%%%%%%%+\frac{1}{\sqrt{k}}\,\big(w_\mathrm{t}+w_\mathrm{h,tg}+w_\mathrm{h,nor}\big)
%%%%%%%%%\frac{w_\mathrm{t}}{\sqrt{k}}+\frac{w_\mathrm{h,tg}}{\sqrt{k}}+\frac{w_\mathrm{h,nor}}{\sqrt{k}}
%%%%%%%%%%\right)
%%%%%%%%%%%%%%%%%%%%%\nonumber
\end{eqnarray}
where $v=v_\mathrm{t}+v_\mathrm{h,nor}\in T_{\mu_g\circ\pi\circ \widetilde{\sigma} (x'_0)}M$;
the latter is unitarily identified with $\mathbb{C}^\mathrm{d}$ by means
of the given Heisenberg local coordinates. Taylor expanding,
%at $y=\mu_g\circ\widetilde{\sigma}_{lj}
%(x'_0)+\frac{v_\mathrm{t}}{\sqrt{k}}+\frac{v_\mathrm{h,nor}}{\sqrt{k}}$,
we obtain with $m_0'=\pi_0(q_0')\in V_{lj}\subseteq F_l$:
\begin{eqnarray}
\label{eqn:taylor-expansion}
f(y)=f\circ \pi(y)&\sim& f\big(\mu_g\circ\sigma (m'_0)\big)+\sum_{j\ge 1}k^{-j/2}\,f_j(v,w),\\
\tau_{lj}(y)=\tau_{lj}\big(\pi(y)\big)&\sim &\tau_{lj}^0\big(m'_0\big)+\sum_{j\ge 1}k^{-j/2}\,\tau_j(v),\nonumber\\
\varsigma_k(y)=\varsigma_k\big(\pi(y)\big)&\equiv&1\,\,\,\,\,\,\,\,\,
\mathrm{if}\,\,
\|v_\mathrm{t}+v_\mathrm{h,nor}\|\lesssim k^{1/10}.\nonumber
\end{eqnarray}
Let us now recall the asymptotic expansion for scaling limits of equivariant Szeg\"{o} kernels proved in
\cite{p}. Given $g_0\in G$, $e^{i\vartheta_0}\in S^1$, $x\in \left (\Phi\circ \pi\right)^{-1}(0)$ and $w=w_\mathrm{t}+w_\mathrm{v}+w_\mathrm{h},
v=v_\mathrm{t}+v_\mathrm{v}+v_\mathrm{h}\in T_{\pi(x)}(M)$,
as $k\rightarrow +\infty$ we have
an asymptotic expansion
\begin{eqnarray}%\label{eqn:expansion-for-pi-equivariant}
\label{eqn:asympt-expansion-general-case}
\lefteqn{\Pi_{\varpi,k}\left
(\mu_{g_0}\circ r_{e^{i\vartheta_0}}\left(x+\frac{w}{\sqrt{k}}\right),x+\frac{v}{\sqrt{k}}\right)}\\
 &\sim& \left
(\frac{k}{\pi}\right
)^{\mathrm{d}-\mathrm{g}/2}\,A_{\varpi,k}(x,g_0,h_0)\,e^{Q(w_\mathrm{v}+w_\mathrm{t},
v_{\mathrm{v}}+v_\mathrm{t})}\,e^{\psi
_2(w_\mathrm{h},v_\mathrm{h})}\cdot \left (1+\sum _{j\ge 1}a_{\varpi
j}(x,w,v)\,k^{-j/2}\right),\nonumber\end{eqnarray}
where in the present situation
\begin{itemize}
  \item $
A_{\omega,k}(x,g_0,h_0)=:2^{\mathrm{g}/2}\,\dim (V_\varpi)\,V_{\mathrm{eff}}(x)^{-1}
\,
\chi _\varpi \left(g_0^{-1}\right)\,
e^{ik\vartheta_0}$,
where $V_{\mathrm{eff}}(x)=V_{\mathrm{eff}}\big(\pi(x)\big)$ is the \textit{effective volume} of $\pi(x)$, that is,
the volume of the $G$-orbit $G\cdot \pi(x)\subseteq M$ with respect to the Riemannian density \cite{bg};
being $G$-invariant, $V_{\mathrm{eff}}$ descends to a smooth function on $M_0$,
again denoted $V_{\mathrm{eff}}$.
\item $
Q(w_\mathrm{v}+w_\mathrm{t},
v_{\mathrm{v}}+v_\mathrm{t})=:-\|v_\mathrm{t}\|^2-\|w_\mathrm{t}\|^2+i\,\big[\omega
_m(w_\mathrm{v},w_\mathrm{t})-\omega
_m(v_\mathrm{v},v_\mathrm{t})\big]$, where $\omega =\frac i2\,\Theta$,
and the norms are taken in the
Hermitian structure $h$ of $TM$ induced by $\omega$.
\item $\psi _2\big(w_\mathrm{h},v_\mathrm{h}\big)=:h_m\big(w_\mathrm{h},v_\mathrm{h}\big)-
  \frac 12\,\left (\|w_\mathrm{h}\|^2+\|v_\mathrm{h}\|^2\right
)$.
\item the $a_{\varpi
j}$'s are polynomials in $v$, $w$ whose coefficients depend on $x$ and $\varpi$.
\item Let $R_N(x,w,v)$ be the remainder term following the first $N$ summands in (\ref{eqn:asympt-expansion-general-case}); then for
$\|w\|,\|v\|\lesssim k^{1/6}$ we have the \lq large ball estimate\rq:
\begin{equation}
\label{eqn:large-ball-estimate}
\big |R_N(x,w,v)\big |\le C_N\,k^{\mathrm{d}-(\mathrm{g}+N+1)/2}\,e^{-\frac{1-\epsilon}{2}
\big(\|w_\mathrm{h}-v_\mathrm{h}\|^2+2\|v_\mathrm{t}\|^2+2\|w_\mathrm{t}\|^2\big)}.
\end{equation}
\end{itemize}

To obtain an asymptotic expansion for
$\Pi_{\varpi,k}\left (\gamma_X^{-1}(y),y\right)$ in
(\ref{eqn:trace-lj-th}), with $y$ as in
(\ref{eqn:x-e-y-riscalati}), we need the the Heisenberg
local coordinates of $\gamma_X^{-1}(y)$. As a first step,
let us work out the underlying preferred coordinates of
$\gamma^{-1}\big(\pi(y)\big)$.

More precisely, recall that our construction involves a moving
Heisenberg local chart $\Psi=\Psi_{lj}$ as in
(\ref{eqn:heisenberg-mobile}), now with $(V_{lj},\sigma_{lj})$ in
place of $(V,\sigma)$. Underlying $\Psi$, there is a moving
preferred local chart
\begin{equation}
\label{eqn:preferred-mobile}
\mathfrak{F}:V_{lj}\times B_{2\mathrm{d}}(\epsilon)\rightarrow M,
\end{equation}
such that for any $m_0'\in V_{lj}$ the partial map
$\mathfrak{f}^{(m_0')}=:\mathfrak{F}(m_0',\cdot,\cdot)$ is a preferred local
chart for $M$ centered at $\sigma(m_0')$, with image containing some
fixed open subset $U\subseteq M$. It will simplify our exposition to
assume, as we may, that for every $m_0'\in V$ we have
% the exponential map
$$\mathfrak{f}^{(m_0')}=\exp_{\sigma(m_0')}\circ \varrho_{m_0'}:B_{2d}(\epsilon)\subseteq\mathbb{C}^\mathrm{d}\cong
T_{\sigma(m_0')}M\rightarrow M,$$ for an appropriate $\epsilon>0$;
here $\varrho_{(m_0')}:\mathbb{C}^\mathrm{d}\cong T_{\sigma(m_0')}M$
is a smoothly varying unitary isomorphism, induced by the choice of
an orthonormal frame for $T^{(1,0)}M$ on a neighborhood of
$\sigma_{lj}(V_{lj})$.

The same then holds for every pair $(g,m_0')\in G\times V_{lj}$.
More precisely, since $\mu_g:M\rightarrow M$ is a Riemannian
isometry, for every $(g,m_0')\in G\times V_{lj}$ the composition
$$\varrho_{(g,m_0')}=:d_{\sigma(m_0')}\mu_g\circ\varrho_{m_0'}:\mathbb{C}^\mathrm{d}\rightarrow
T_{\sigma(m_0')}M$$ is unitary, where write $\sigma=\sigma_{lj}$,
and $\mu_g\circ \exp_{\sigma(m_0')}=\exp_{\mu_g\circ
\sigma(m_0')}\circ d_{\sigma(m_0')}\mu_g$. By construction,
\begin{eqnarray*}
\lefteqn{\mathfrak{f}^{(g,m_0')}=\mu_g\circ \mathfrak{f}^{(m_0')}=\mu_g\circ
\exp_{\sigma(m_0')}\circ
\varrho_{(m_0')}}\\
&=&\exp_{\mu_g\circ\sigma_{lj}(m_0')}\circ
d_{\sigma(m_0')}\mu_g\circ
\varrho_{(m_0')}=\exp_{\mu_g\circ\sigma(m_0')}\circ
\varrho_{(g,m_0')}.
\end{eqnarray*}
Given $(g,m_0')\in G\times V_{lj}$ and $v\in T_{\mu_g\circ
\sigma(m_0')}M$ with $\|v\|<\epsilon$, we shall set $\mu_g\circ
\sigma_{lj}(m_0')+v=:\mathfrak{f}^{(g,m_0')}\circ
\varrho_{(g,m_0')}^{-1}(v)$.

For every $m\in \widetilde{\mathrm{Fix}(\gamma_0)}$ there exists a
unique $g_m\in G$ such that $\gamma(m)=\mu_{g_m}(m)$ (Definition \ref{defn:componenti-luogo-fisso}).
Hence $\forall\,(g,m_0')\in G\times M$, with $\kappa=:g\,g_{\sigma
(m_0')}^{-1}\,g^{-1}$,
\begin{equation}
\label{eqn:gamma-lungo-Fl}
\gamma^{-1}\Big( \mu_g\circ \sigma (m_0')\Big)=\mu_{\kappa}\circ
\mu_g\circ \sigma (m_0')=\mu_{\kappa g}\circ
 \sigma (m_0').
\end{equation}
Since $\gamma:M\rightarrow M$ is also a Riemannian isometry, by (\ref{eqn:gamma-lungo-Fl})
we also have
\begin{eqnarray}
\label{eqn:gamma-isometry} \lefteqn{\gamma^{-1}\circ
\exp_{\mu_g\circ \sigma (m_0')}}\\
&=&\exp_{\gamma^{-1}\circ\mu_g\circ \sigma (m_0')}\circ
d_{\mu_g\circ \sigma (m_0')}\gamma^{-1}=\exp_{\mu_{\kappa g}\circ \sigma
(m_0')}\circ d_{\mu_g\circ \sigma (m_0')}\gamma^{-1}.\nonumber
\end{eqnarray}
With the previous convention, (\ref{eqn:gamma-isometry}) implies
that $\forall\,v\in T_{\mu_g\circ \sigma (m_0')}M$
\begin{equation}
\label{eqn:azione-di-gamma} \gamma^{-1}\Big( \mu_g\circ \sigma
(m_0')+v\Big)=\mu_{\kappa g}\circ\sigma (m_0')+d_{ \mu_g\circ \sigma
(m_0')}\gamma^{-1}(v).
\end{equation}
Lifting this to $X$, we deduce that $\forall\,
\left(g,x_0'\right)\in G\times X(V_{lj})$
\begin{equation}
\label{eqn:gamma-X-eq-locali}
\gamma_X^{-1}\Big( \mu_g\circ
\widetilde{\sigma}
(x_0')+v\Big)=r_{e^{i\beta(g,x_0',v)}}\Big(\mu_{\kappa g}\circ\widetilde{\sigma}
(x_0')+d_{ \mu_g\circ \sigma (m_0')}\gamma^{-1}(v)\Big),
\end{equation}
with $m_0'=:\pi_0(x_0')\in V_{lj}\subseteq F_l$, for an appropriate smooth real function
$\beta:G\times X(V_{lj})\times B_{2\mathrm{d}}(\epsilon)\rightarrow \mathbb{R}$,
uniquely determined up to an integer multiple of $2\pi$. To determine
$\beta$, recall that $h_l=e^{i\theta_l}\in S^1$ is uniquely determined by the condition
$\widetilde{\gamma}_0\big((r,a)\big)=(r,h_l\,a)$, $\forall\,r\in F_l\subseteq M_0,\,(r,l)\in A_0(r)$;
here $A_0(r)$ is the fiber of $A_0$ at $r$, and $\widetilde{\gamma}_0:A_0\rightarrow A_0$
is the linearization of $\gamma_0$ (Definition \ref{defn:quotient-linearization}).

\begin{lem}
\label{lem:local-behav-theta}
Perhaps after adding a suitable integer multiple of $2\pi$,
we may assume that $\beta-\theta_l$ vanishes to third order at $v=0$, that is,
$$
\beta(g,x_0',v)=\theta_l+\sum _{|I|+|J|=3}c_{I,J}(g,x_0')v^I\,\overline{v}^J+R(g,x_0',v),
$$
where $R(g,x_0',\cdot)$ vanishes to fourth order at $v=0$.
\end{lem}

\textit{Proof.}
Let $\gamma_{X_0}:X_0\rightarrow X_0$ be the contactomorphism induced by $\gamma_X$ by passage
to the quotient. In other words, $\gamma_{X_0}$ is the restriction to $X_0$ of the dual linearization
$\left(\widetilde{\gamma}_0^{-1}\right)^t$ on $A^*_0$. Let us momentarily write
$x_0'=(m_0',\eta)$, where $m_0'=\pi(x_0')\in F_l\subseteq \mathrm{Fix}(\gamma_0)$, and
$\eta\in A_0^*(m_0')$ has unit norm. We obtain
\begin{eqnarray}
\label{eqn:linarizzazione-duale-quoziente}
 \gamma_{X_0}^{-1}(x_0')=\left (\gamma_0^{-1}(m_0'),\eta\circ \widetilde{\gamma}_0\right)
 =\left (m_0',e^{i\theta_l}\eta\right)=r_{e^{i\theta_l}}(x_0').
 \end{eqnarray}
On the other hand, (\ref{eqn:gamma-X-eq-locali}) with $v=0$ descends on $X_0$ to the relation
\begin{equation}\label{eqn:relazione-indotta-quoziente}
   \gamma_{X_0}^{-1}
(x_0')=r_{e^{i\beta(g,x_0',0)}}(x_0').
\end{equation}
Now (\ref{eqn:linarizzazione-duale-quoziente}) and (\ref{eqn:relazione-indotta-quoziente})
imply that $\beta(g,x_0',0)-\theta_l=2\pi c$ for some $c\in \mathbb{N}$; by continuity, $c$ is
constant, and we may assume without loss that $c=0$.

Next, we make use of the fact that $\gamma_X$ is a contactomorphism, that is, $\gamma_X^*(\alpha)=\alpha$.

Let the Heisenberg local chart $\psi^{(g,x_0')}$ be as in (\ref{eqn:heisenber-local-chart-0}),
with image an open neighborhood $X^{(g,x_0')}\subseteq X$ of $\mu_g\circ \widetilde{\sigma}(x_0')$.
Let
$\left(\mathbf{z}^{(g,x_0')},\vartheta^{(g,x_0')}\right):X^{(g,x_0')}\rightarrow \mathbb{C}^\mathrm{d}\times
\mathbb{R}$ be the
%=\left (z_1^{(g,x_0')},\ldots,z_\mathrm{d}^{(g,x_0')},\vartheta^{(g,x_0')}\right)
associated local coordinates, as in
(\ref{eqn:heisenber-local-chart}).
We write $\mathbf{z}^{(g,x_0')}=\mathbf{a}^{(g,x_0')}+i\mathbf{b}^{(g,x_0')}$,
with $\mathbf{a}^{(g,x_0')},\,\mathbf{b}^{(g,x_0')}$ real-valued. Then by \cite{sz}, \S 1 the local expression
of $\alpha$ on $X^{(g,x_0')}$ has the form
$$
\alpha=d\vartheta^{(g,x_0')}+\mathbf{a}^{(g,x_0')}\,d\mathbf{b}^{(g,x_0')}-\mathbf{b}^{(g,x_0')}\,d\mathbf{a}^{(g,x_0')}
+\beta^{(g,x_0')}\left (\mathbf{z}^{(g,x_0')}\right),
$$
where $\beta^{(g,x_0')}\left (\mathbf{z}^{(g,x_0')}\right)=O\left (\left\|\mathbf{z}^{(g,x_0')}\right\|^2\right)$.
A similar expression, with $g$ replaced by $\kappa\,g$, holds for $\alpha$ on $X^{(\kappa g,x_0')}$. Since
$\left(\gamma_X^{-1}\right)^*\left(\vartheta^{(\kappa g,x_0')}\right)=\vartheta^{(g,x_0')}+\beta$, and
$d_{ \mu_g\circ \sigma (m_0')}\gamma^{-1}$ is unitary, hence symplectic,
(\ref{eqn:gamma-X-eq-locali}) implies
\begin{eqnarray*}
\lefteqn{\alpha=d\vartheta^{(g,x_0')}+\mathbf{a}^{(g,x_0')}\,d\mathbf{b}^{(g,x_0')}-\mathbf{b}^{(g,x_0')}\,d\mathbf{a}^{(g,x_0')}
+\beta^{(g,x_0')}\left (\mathbf{z}^{(g,x_0')}\right)}\\
&=&\left(\gamma_X^{-1}\right)^*\left ( d\vartheta^{(\kappa g,x_0')}+\mathbf{a}^{(\kappa g,x_0')}\,d\mathbf{b}^{(\kappa g,x_0')}-\mathbf{b}^{(\kappa g,x_0')}\,d\mathbf{a}^{(\kappa g,x_0')}
+\beta^{(\kappa g,x_0')}\left (\mathbf{z}^{(\kappa g,x_0')}
\right)\right)\\
&=&d\vartheta^{(g,x_0')}+d\beta+\mathbf{a}^{(g,x_0')}\,d\mathbf{b}^{(g,x_0')}-\mathbf{b}^{(g,x_0')}\,d\mathbf{a}^{(g,x_0')}
+\left(\gamma_X^{-1}\right)^*\left (\beta^{(\kappa g,x_0')}\left (\mathbf{z}^{(\kappa g,x_0')}
\right)\right).
\end{eqnarray*}
We deduce
$$
d\beta=\beta^{(g,x_0')}\left (\mathbf{z}^{(g,x_0')}\right)-\left(\gamma_X^{-1}\right)^*\left (\beta^{(kg,x_0')}\left (\mathbf{z}^{(kg,x_0')}
\right)\right)=O\left (\left\|\mathbf{z}^{(g,x_0')}\right\|^2\right).
$$
\hfill Q.E.D.

\begin{cor}
\label{eqn:asympt-of-the-angle}
Let $y=y_k$ ($k=1,2,\ldots$) be as in (\ref{eqn:x-e-y-riscalati}), and
set $\kappa=:g\,g_{\sigma
(m_0')}^{-1}\,g^{-1}\in G$, $v=:v_\mathrm{t}+v_{\mathrm{h,nor}}\in T_{\mu_g\circ \sigma (m_0')}M$.
In the Heisenberg local chart $\psi^{(\kappa g,x_0')}$,
$$
\gamma_X^{-1}(y)=r_{e^{i\beta_k(g,x_0',v)}}\left(\mu_{\kappa g}\circ\widetilde{\sigma}
(x_0')+\frac{1}{\sqrt{k}}\,d_{ \mu_g\circ \sigma (m_0')}\gamma^{-1}(v)\right),
$$
where
$\beta_k(g,x_0',v)\sim \theta_l+\sum _{j\ge 0}k^{-(3+j)/2}\,b_j(g,x_0',v)$ as $k\rightarrow +\infty$.
\end{cor}

On the upshot, with $y$ as in (\ref{eqn:x-e-y-riscalati}), we obtain
\begin{eqnarray}
\label{eqn:asymptotic-expansion-2nd-factor}
\lefteqn{\Pi_{\varpi,k}\left (\gamma_X^{-1}(y),y\right)}\\
&=&\Pi_{\varpi,k}\left (r_{e^{i\beta_k(g,x_0',v)}}\left(\mu_{\kappa g}\circ\widetilde{\sigma}
(x_0')+\frac{d_{ \mu_g\circ \sigma (m_0')}\gamma^{-1}(v)}{\sqrt{k}}\right),\mu_g\circ\widetilde{\sigma} (x'_0)
+\frac{v}{\sqrt{k}}\right)\nonumber\\
&=&e^{ik\beta_k(g,x_0',v)}\,\Pi_{\varpi,k}\left (\mu_{\kappa g}\circ\widetilde{\sigma}
(x_0')+\frac{d_{ \mu_g\circ \sigma (m_0')}\gamma^{-1}(v)}{\sqrt{k}},\mu_g\circ\widetilde{\sigma} (x'_0)
+\frac{v}{\sqrt{k}}\right),\nonumber\end{eqnarray}
where $v=v_\mathrm{t}+v_{\mathrm{h,nor}}\in T_{\mu_g\circ \sigma (m_0')}M$.
Now
$$
\mu_g\circ\widetilde{\sigma} (x'_0)
+\frac{v}{\sqrt{k}}
%=\mu_{h\kappa^{-1}}\circ \mu_\kappa\left (\mu_g\circ\widetilde{\sigma} (x'_0)
%+\frac{w}{\sqrt{k}}\right)
=\mu_{\kappa ^{-1}}\left (\mu_{\kappa g}\circ\widetilde{\sigma} (x'_0)
+\frac{1}{\sqrt{k}}\,d_{\mu_g\circ \sigma(m_0')}\mu_\kappa (v)\right)
$$
Now we remark that $v_\mathrm{v}=v_{\mathrm{h,tg}}=0$; furthermore, $d\mu_\kappa$ and
$d\gamma^{-1}$ preserve the decomposition of tangent vectors described in (\ref{eqn:struttura-locale-M-phi}) - (\ref{eqn:struttura-locale-M-phi-vert-1}),
and the norm of each component.
Let $v_0\in T_{m_0'}M_0$
correspond to
$v_\mathrm{h}\in H_{\mu_g\circ \sigma (m_0')}\subseteq T_{\mu_g\circ \sigma (m_0')}M$.
Any $r\in T_{m_0'}M_0$ may be decomposed as
$r=r_{\mathrm{tg}}+r_{\mathrm{nor}}$, where
$r_{\mathrm{tg}}\in T_{m_0'}F_l$, $r_{\mathrm{nor}}\in \big(T_{m_0'}F_l\big)^\perp$.
In our case, $v_0=v_{0,\mathrm{nor}}$, whence
$$\left (d_{ \mu_g\circ \sigma (m_0')}\gamma^{-1}(v)\right)_0=d_{m_0'}\gamma_0^{-1}(v_{0,\mathrm{nor}})=
d_{m_0'}\gamma_0^{-1}\big(v_{0}\big)_\mathrm{nor}.$$
Thus,
\begin{eqnarray}
\label{eqn:psi-2-sotto}
%%%%\psi_2\big(w_\mathrm{h},v_\mathrm{h}\big)&=&-\frac 12\,\big\|w_{0,\mathrm{tg}}\big\|^2+\psi_2\big(w_{0,\mathrm{nor}},
%%%%%v_{0,\mathrm{nor}}),\\
\psi_2\Big(d_{\mu_g\circ \sigma(m_0')}\mu_\kappa (v)_\mathrm{h},d_{ \mu_g\circ \sigma (m_0')}\gamma^{-1}(v)_\mathrm{h}\Big)
=
\psi_2\Big(v_{0,\mathrm{nor}},d_{m_0'}\gamma_0^{-1}(v_{0,\mathrm{nor}})\Big).
\end{eqnarray}
Therefore, by (\ref{eqn:asympt-expansion-general-case}) we deduce
\begin{eqnarray}
\label{eqn:asympt-exp-second-factor}
\lefteqn{\Pi_{\varpi,k}\left (\mu_{\kappa ^{-1}}\left (\mu_{\kappa g}\circ\widetilde{\sigma} (x'_0)
+\frac{1}{\sqrt{k}}\,d_{\mu_g\circ \sigma(m_0')}\mu_\kappa (v)\right),\mu_{\kappa g}\circ\widetilde{\sigma}
(x_0')+\frac{1}{\sqrt{k}}\,d_{ \mu_g\circ \sigma (m_0')}\gamma^{-1}(v)\right)}\nonumber\\
&\sim& \left
(\frac{k}{\pi}\right
)^{\mathrm{d}-\mathrm{g}/2}\,2^{\mathrm{g}/2}\,\frac{\dim (V_\varpi)}{V_{\mathrm{eff}}(x_0')}
\,
\chi _\varpi \left(\kappa \right)\,e^{-2
\|v_\mathrm{t}\|^2}\,e^{\psi_2\big(v_{0,\mathrm{nor}},d_{m_0'}\gamma_0^{-1}(v_{0,\mathrm{nor}})\big)}\nonumber\\
&&\cdot \left (1+\sum _{j\ge 1}a_{\varpi
j}(x,v)\,k^{-j/2}\right).
\end{eqnarray}
Taking conjugates, we obtain from (\ref{eqn:asymptotic-expansion-2nd-factor}) and (\ref{eqn:asympt-exp-second-factor})
that for $y$ given by (\ref{eqn:x-e-y-riscalati}) we have
\begin{eqnarray}
\label{eqn:asymptotic-expansion-2nd-factor-final}
\lefteqn{\Pi_{\varpi,k}\left (\gamma_X^{-1}(y),y\right)} \\
&\sim& \left
(\frac{k}{\pi}\right
)^{\mathrm{d}-\mathrm{g}/2}\,2^{\mathrm{g}/2}\,\frac{\dim (V_\varpi)}{V_{\mathrm{eff}}(x_0')}
\,e^{ik\theta_l}\,
\chi _\varpi \left(g_{\sigma(m_0')}\right)\,e^{-2
\|v_\mathrm{t}\|^2}\,e^{\psi_2\big(d_{m_0'}\gamma_0^{-1}(v_{0,\mathrm{nor}}),
v_{0,\mathrm{nor}}\big)}\nonumber\\
&&\cdot \left (1+\sum _{j\ge 1}b_{\varpi
j}(x,v)\,k^{-j/2}\right).   \nonumber
\end{eqnarray}

Given (\ref{eqn:taylor-expansion}), (\ref{eqn:asymptotic-expansion-2nd-factor-final}),
and Corollary \ref{eqn:asympt-of-the-angle},
as $k\rightarrow +\infty$ the integrand of (\ref{eqn:trace-lj-th}) (omitting $\varsigma_k(x,y)$)
admits an asymptotic expansion
\begin{eqnarray}
\label{eqn:asympt-exp-integrand}
\lefteqn{\tau_{lj}(y)\,\Pi_{\varpi,k}\left (\gamma_X^{-1}(y),y\right)\,f(y)}\\
&\sim& \tau_{lj}^0(m_0')\,f\big(\mu_g\circ \sigma (m'_0)\big)\,\left
(\frac{k}{\pi}\right
)^{\mathrm{d}-\mathrm{g}/2}\,2^{\mathrm{g}/2}\,\frac{\dim (V_\varpi)}{V_{\mathrm{eff}}(x_0')}
\,e^{ik\theta_l}
\chi _\varpi (F_l)\,e^{-
2\|v_\mathrm{t}\|^2}\nonumber\\
&&\cdot e^{\psi_2\big(d_{m_0'}\gamma_0^{-1}(v_{0,\mathrm{nor}}),
v_{0,\mathrm{nor}}\big)}\,\left (1+\sum _{j\ge 1}c_{\varpi
j}(x,v)\,k^{-j/2}\right),   \nonumber
\end{eqnarray}
for certain polynomials $c_{\varpi
j}$ in $v$; here $\chi _\varpi (F_l)$ is as in Definition \ref{defn:componenti-luogo-fisso}.

(\ref{eqn:x-e-y-riscalati}) is obtained by composing the parametrization $\mathcal{T}=\mathcal{T\mathcal{}}_{lj}$
in (\ref{eqn:parametr-per-tilde-b}) with
a rescaling by $k^{-1/2}$. Let us denote  this composition by $\mathcal{T}_{lj}^{(k)}$.
Accordingly, (\ref{eqn:trace-lj-th}) may be viewed as an integral over
$G\times X(V)\times \mathbb{R}^\mathrm{g}\times \mathbb{C}^{\mathrm{c}_l}$, with respect to
the density $\left(\mathcal{T}_{lj}^{(k)}\right)^*\big(\mathrm{dens}_{X}\big)$.
In view of Lemma 3.9 of \cite{dp}, and by the construction of $\mathcal{T}_{lj}^{(k)}$ using Heisenberg coordinates
and rescaling, this admits an asymptotic expansion of the form
\begin{eqnarray}
\label{eqn:asympt-exp-density}
\lefteqn{\left(\mathcal{T}_{lj}^{(k)}\right)^*\big(\mathrm{dens}_{X}\big)(y)}\\
&\sim&k^{-\mathrm{c}_l-\mathrm{g}/2}\,V_{\mathrm{eff}}\left(m_0'\right)
\, \big|d\nu(g)\big|\,\mathrm{dens}_{X_0}(x_0')\,\big| dv\big|
\cdot \left (1+\sum _{a\ge 1}k^{-a/2}d_a(\theta,x_0',v)\right),\nonumber
\end{eqnarray}
where $\big| dv\big|$ is the Lebesgue measure on $\mathbb{R}^\mathrm{g}\times \mathbb{C}^{\mathrm{c}_l}$,
and each $d_a(\theta,x_0',v)$ is a polynomial in $v$.

Before proceeding, we need to establish that the asymptotic expansion obtained by multiplying
(\ref{eqn:asympt-exp-integrand}) and
(\ref{eqn:asympt-exp-density}) can be integrated
term by term; to this end, let us pause on the remainder term.
By the considerations preceding (\ref{eqn:x-e-y-riscalati}), integration $v$ is over a ball
of radius $\approx k^{1/6}$ in
$\mathbb{R}^{\mathrm{g}}\times \mathbb{C}^{\mathrm{c}_l}$. On the domain of integration, therefore,
the remainder term in (\ref{eqn:asymptotic-expansion-2nd-factor-final}) satisfies
(\ref{eqn:large-ball-estimate}), with $w_\mathrm{h}$ replaced by $d\gamma^{-1}(v_\mathrm{h})$.
When we multiply the asymptotic expansions, therefore, one of the typical contributions due to the remainder terms
is bounded by
\begin{eqnarray}
\label{eqn:first-contribution-remainder}
C\left (\frac k \pi\right)^{\mathrm{d}-g/2-(N+1)/2}\,p_N(v)\,e^{-2(1-\epsilon)\|v_\mathrm{t}\|^2
-\frac{1-\epsilon}{2}\,\|v_\mathrm{h,nor}-d\gamma^{-1}(v_\mathrm{h,nor})\|^2},
\end{eqnarray}
where $p_N$ is some polynomial in $v=v_\mathrm{h}+v_\mathrm{t}$, and $N$ is a
positive integer, that may be assumed to grow to infinity with the length of our expansion.
On the other hand,
\begin{eqnarray}
\label{eqn:second-string-of-estimates}
\left\|v_\mathrm{h,nor}-d\gamma^{-1}(v_\mathrm{h,nor})\right\|\ge D\,\|v_\mathrm{h,nor}\|,
\nonumber
\end{eqnarray}
where $D^{-1}>0$ is the operator norm of $\left(I-d\gamma^{-1}\right)^{-1^{}}$
acting on any fiber of the normal bundle $N_l$.
The other terms can be handled
in a similar way.

We can thus integrate term by term, and this proves the existence of
an asymptotic expansion for $\mathrm{trace}\left (\Psi_{\varpi,k}\right)_{lj}$
in (\ref{eqn:trace-lj-th}) as $k\rightarrow +\infty$, and therefore for $\mathrm{trace}\left (\Psi_{\varpi,k}\right)$.
Let us now explicitly compute the leading term.

Since the region where $\zeta_k\neq 1$ yields a contribution to the integral which is
$O\left (k^{-\infty}\right)$, in the following we shall set $\zeta_k=1$.
%Let us first integrate in the transverse directions, $w_{\mathrm{t}}$ and $v_{\mathrm{t}}$, and in the horizontal.

We have%
\begin{equation}
\label{eqn:integr-in-w-transv}
\int _{\mathbb{R}^{g}}\,e^{-2\|v_{\mathrm{t}}\|^2}\,
dv_{\mathrm{t}}=\int _{\mathbb{R}^{g}}e^{-2\|\mathbf{x}\|^2}\,d\mathbf{x}=\left (
\frac \pi 2\right)^{\mathrm{g}/2}.
\end{equation}

Next, let $\Lambda_l\in U(\mathrm{c}_l)$ denotes the unitary matrix representing
the restriction of $d_{m_0'}\gamma_0$ to the normal space
of $F_l$ at $m_0'$, in the induced coordinates. The conjugacy class of $\Lambda_l$
only depends on $l$.
Let $\big(v_1,\ldots,v_{\mathrm{c}_l}\big)$ be an orthonormal basis of $\mathbb{C}^{\mathrm{c}_l}$ composed of eigenvectors
of $\Lambda_l$, with corresponding eigenvalues $\lambda_1,\ldots,
\lambda_{\mathrm{c}_l}\in S^1\setminus \{1\}$.
If $v_{0,\mathrm{nor}}=\sum _{j=1}^{\mathrm{c}_l}a_j\,v_j$ we have
\begin{equation}\label{eqn:psi-2-sviluppo}
\psi_2\big(\Lambda_l^{-1} v_{0,\mathrm{nor}},
v_{0,\mathrm{nor}}\big)
    =\sum _{j=1}^{\mathrm{c}_l}\left(\overline{\lambda}_j-1\right)\,|a_j|^2.
\end{equation}
Therefore, recalling that $\int _\mathbb{C}e^{s|u|^2}\,du=-\pi/s$ if $\Re(s)<0$, we
get

\begin{eqnarray}
\label{eqn:third-factor-szego-second-int}
\lefteqn{
\int _{\mathbb{C}^{\mathrm{c}_l}}
e^{
\psi_2\big(d_{m_0'}\gamma_0^{-1}(v_{0,\mathrm{nor}}),
v_{0,\mathrm{nor}}\big)}\,dv_{0,\mathrm{nor}}
}\\
&=&\int _{\mathbb{C}^{\mathrm{c}_l}}
e^{\sum _{j=1}^{\mathrm{c}_l}\left(\overline{\lambda}_j-1\right)\,|a_j|^2
}\,da=\prod_{j=1}^{\mathrm{c}_l}\int _{\mathbb{C}}
e^{\left(\overline{\lambda}_j-1\right)\,|u|^2
}\,du=\pi^{\mathrm{c}_l}\,\prod_{j=1}^{\mathrm{c}_l}\frac{1}{1-\overline{\lambda}_j}
\nonumber\\
&=&\frac{\pi^{\mathrm{c}_l}}{\det\big(
\mathrm{id}_{\big(N_l\big)_{m_0'}}-
\left.d_{m_0'}\gamma_0^{-1}\right|_{N_{l,m_0'}}
\big)
}=\frac{\pi^{\mathrm{c}_l}}{c_l(\gamma)},\nonumber
\end{eqnarray}
where $c_l(\gamma)$ is as in (\ref{eqn:determinant-factor}).
On the upshot, the leading term of the asymptotic expansion for (\ref{eqn:trace-lj-th}) is:
\begin{eqnarray}\label{eqn:trace-lj-th-asymptotics}
\left (\frac{k}{\pi}\right)^{\mathrm{d}_l}\,\frac{e^{ik\theta_l}}{c_l(\gamma)}\,
\dim(V_\varpi)\,\chi_\varpi(F_l)\cdot\int _{X_0(V_{lj})}
\tau_{lj}(m_0')\cdot\left(\int _{G}f\Big(\mu_g\circ \sigma (m'_0)\Big)\,d\nu(g)\right)\,\mathrm{dens}_{X_0}(x_0')\nonumber\end{eqnarray}
where $m_0'=\pi_0(x_0')\in V_{lj}$.
To complete the proof of Theorem \ref{thm:main} we need only sum over $l,j$.

\hfill Q.E.D.

\end{document}